\documentstyle[12pt]{article}
\newlength{\abstractwidth}
\renewcommand{\thefootnote}{\fnsymbol{footnote}}
\renewcommand{\thanks}[1]{\footnote{#1}} % Use this for footnotes
\newcommand{\starttext}{ \setcounter{footnote}{0}
\renewcommand{\thefootnote}{\arabic{footnote}}}

\newcommand{\be}{\begin{equation}}
\newcommand{\bea}{\begin{eqnarray}}
\newcommand{\eea}{\end{eqnarray}} \newcommand{\ee}{\end{equation}}
\newcommand{\N}{{\cal N}} \newcommand{\<}{\langle}
\renewcommand{\>}{\rangle} \def\ba{\begin{eqnarray}}
\def\ea{\end{eqnarray}}

\def\D{{\cal D}}

\def\N{{\cal N}}

\def\o{\omega}
\def\Re{{\rm Re}}

\def\log{\,{\rm log}\,}

\def\o{\omega}

\def\b{\beta}

\def\e{\varepsilon}

\def\o{\omega}

\def\D{\Delta}

\def\N{\bf N}

\def\R{{\bf R}}

\def\p{\partial}

\def\ddb{{\partial\bar\partial}}

\def\D{\Delta}

\def\[{{\bf [}}
\def\]{{\bf ]}}

\begin{document}
\starttext \baselineskip=18pt \setcounter{footnote}{0}
\newtheorem{theorem}{Theorem}
\newtheorem{lemma}{Lemma}
\newtheorem{corollary}{Corollary}
\newtheorem{definition}{Definition}
\newtheorem{conjecture}{Conjecture}
\newtheorem{proposition}{Proposition}

\begin{center}
{\Large \bf ON ESTIMATES FOR THE FU-YAU GENERALIZATION OF A STROMINGER SYSTEM
\footnote{Work supported in part by the National Science Foundation under Grant DMS-12-66033 and DMS-1308136. 
Keywords: Hessian equations; symmetric functions of eigenvalues; Moser iteration,
inequalities of Guan-Ren-Wang; maximum principles. AMS classification numbers: 32Q26 (32Q15, 32Q20, 32U05, 32W20), 35Kxx.}}

\bigskip

{\large Duong H. Phong, Sebastien Picard, and Xiangwen Zhang} \\

\medskip

\begin{abstract}

\medskip
\small{
We study an equation proposed by Fu and Yau as a natural $n$-dimensional generalization of a Strominger system that they solved in dimension $2$. It is a complex Hessian equation with right hand side depending on gradients. Building on the methods of Fu and Yau, we obtain $C^0$, $C^2$, and $C^{2,\alpha}$ a priori estimates. We also identify difficulties in extending the Fu-Yau arguments for non-degeneracy from dimension $2$ to higher dimensions.
}

\end{abstract}

\end{center}

\baselineskip=15pt
\setcounter{equation}{0}
\setcounter{footnote}{0}

\section{Introduction}
\setcounter{equation}{0}

In 1985, Strominger \cite{S} proposed a system of equations for compactifications of superstring theories which satisfy the key physical requirement of $N=1$ supersymmetry. These equations are also remarkable from the mathematical standpoint, as they combine in a novel way features of Ricci-flat metrics on Calabi-Yau manifolds together with Hermitian-Einstein metrics on holomorphic vector bundles. Solutions of Strominger systems were indeed obtained perturbatively by Li and Yau \cite{LY} from Ricci-flat and Hermitian-Einstein metrics. However, non-perturbative solutions proved to be daunting, and it was a major breakthrough when Fu and Yau \cite{FY} obtained the first such solution, some twenty years after Strominger's original proposal.

\medskip
The particular Strominger solution obtained by Fu and Yau was a toric fibration over a $K3$ surface. For such manifolds, Fu and Yau succeeded in reducing the Strominger system to the special case in dimension $n=2$ of the following equation, on a compact $n$-dimensional K\"ahler manifold $(X,\omega)$,
\bea
\label{FY}
i \ddb (e^u - \alpha f e^{-u}) \wedge \o^{n-1} + n \alpha  i \ddb u \wedge i \ddb u \wedge \o^{n-2} + \mu {\o^n \over n!} =0,
\eea
where $\alpha > 0$ is a constant, $f \geq 0$ is a smooth function, $\mu$ is a smooth function such that $\int_X \mu = 0$, and the ellipticity condition described further below in (\ref{C0_ellipticity}) is imposed. When $n=2$, it becomes a Monge-Amp\`ere equation, and Fu and Yau \cite{FY} suggested the problem of studying the equation (\ref{FY}) for general dimension $n$.

\medskip
In this paper, we provide a partial answer to the problem raised by Fu and Yau. More specifically, we express the equation (\ref{FY}) in a more standard complex Hessian type equation (see (\ref{scalar_eqn})), and we establish $C^0$, $C^2$, and $C^{2,\alpha}$ a priori estimates for the equation. An upper $C^1$ bound is automatic from the equation. But just as in the case $n=2$ treated by Fu and Yau, the $C^2$ estimate is contingent upon a lower bound for the second symmetric function $\sigma_2(g')$ of the eigenvalues of the unknown Hermitian form $g_{\bar kj}'$ given in (\ref{gprime}). Indeed, this is equivalent to an improved gradient estimate. One of the key innovations of Fu and Yau was a proof of such a lower bound in dimension $n=2$. However, while we were able to obtain a sharp generalization of their computations to arbitrary dimensions, it turned out that this was not strong enough to imply the desired lower bound (see \S 7),
and it is at this time unclear whether such a lower bound does hold.

\medskip
Our proof of the $C^0$ estimate is a close parallel of the proof by Moser iteration methods used in \cite{FY}. The $C^2$ estimate also builds in an essential way on the methods of \cite{FY}, but we also exploit some new inequalities due to Guan, Ren, and Wang \cite{GRW} in their work on real Hessian equations with gradient terms on the right hand side. Although the $C^{2,\alpha}$ estimate does not require much new work, it does not follow from the classical Evans-Krylov theory due to the dependence of the gradient on the right hand side. However, we can obtain the desired estimate by using the recent works of Wang \cite{Wang} and Tosatti-Wang-Weinkove-Yang \cite{TWWY} which deal with the $C^{2, \alpha}$ regularity of complex Monge-Amp\`ere type equations with H\"older regular right hand side. The estimate still open is the lower bound for $\sigma_2(g')$. We discuss in detail the difficulties in trying to extend to higher dimensions
the Fu-Yau arguments for a lower bound for $\sigma_2(g')$. To handle higher dimensions, we work with general coordinate systems rather than the adapted ones with $\nabla u=(u_1,0,\cdots,0)$ used by Fu and Yau. This allows us a simplified and more transparent derivation of the Fu-Yau results for $n=2$, and a clearer picture of why their arguments are not strong enough for higher dimensions. Because of the complexity of the calculations and possibly for future use, this is presented in detail in section \S 7.

\section{The Fu-Yau Equation}
\setcounter{equation}{0}

We begin by writing equation (\ref{FY}) proposed by Fu-Yau \cite{FY} in a more explicit form. Let $\Lambda$ be a Hermitian $(1,1)$-form, and let $\sigma_k(\Lambda)$ be the $k$-th symmetric function of its eigenvalues relative to the K\"ahler form $\omega$, that is,
\bea
\sigma_k(\Lambda)={n\choose k} 
{\Lambda^k\wedge \omega^{n-k}\over\omega^n}
=
\sum_{j_1<\cdots<j_k}\Lambda_{j_1}\cdots\Lambda_{j_k}
\eea
where $\Lambda_j$ denotes the eigenvalues of $\Lambda$ relative to $\omega$. We shall also simplify equation (\ref{FY}) by writing $f$ instead of $\alpha f$ and $\mu$ instead of ${\mu \over (n-2)!}$. Using this notation, equation (\ref{FY}) can then be rewritten as
\bea
\label{FY1}
(n-1)\Delta(e^u-fe^{-u})+
2n\alpha\sigma_2(i\ddb u)+\mu=0.
\eea
The ellipticity condition for this equation is that the
$(1,1)$-Hermitian form $\tilde g_{\bar{j} k}$ defined below be strictly
positive definite,
\bea
\label{C0_ellipticity}
\tilde g_{\bar{j} k}
=
(n-1)(e^u+fe^{-u})g_{\bar{j} k}+
2n\alpha ((\Delta u)g_{\bar{j} k}-u_{\bar{j} k}) >0,
\eea
%Indeed, it is important for geometric applications that $\tilde g_{\bar{j} k}$ be a metric.
where $\o = i \sum g_{\bar j k}dz^k \wedge d\bar z^j$. It is convenient to introduce also the following Hermitian $(1,1)$-form,
\bea\label{gprime}
g'_{\bar{j} k}=(e^u+fe^{-u})g_{\bar{j} k}+2n\alpha u_{\bar{j} k}.
\eea
If we denote by $\lambda_j$ the eigenvalues of $i\ddb u$, by $\lambda_j'$ the eigenvalues of $g_{\bar{j} k}'$, and by $\tilde\lambda_j$ the eigenvalues of $\tilde g_{\bar{j} k}$, all with respect to $g_{\bar{j} k}$, then it is easy to see that
\bea\label{relationg}
\lambda_j'=(e^u+fe^{-u})+2n\alpha\lambda_j,
\qquad
\tilde\lambda_j=\sum_{k\not=j}\lambda_k',
\eea
and hence the following relations between the symmetric functions of $i\ddb u$ and $g_{\bar{j} k}'$,
\bea
\sigma_1(g')&=& n(e^u+fe^{-u})+2n\alpha \sigma_1(i\ddb u)
\nonumber\\
\sigma_2(g')&=&4n^2\alpha^2\sigma_2(i\ddb u)
+
2n(n-1)\alpha (e^u+fe^{-u}) \sigma_1(i\ddb u)
+
{n(n-1)\over 2}(e^u+fe^{-u})^2 \nonumber
\eea
and between the symmetric functions of $g_{\bar{j} k}'$ and $\tilde g_{\bar{j} k}$,
\bea
\sigma_1(\tilde g)&=&(n-1)\sigma_1(g')
\nonumber\\
\sigma_2(\tilde g)&=&{1\over 2}(n-1)(n-2)
\sigma_1(g')^2+\sigma_2(g').
\eea
Substituting equation (\ref{FY1}) in the above expression for $\sigma_2(g')$, we can re-write the equation in terms of $g_{\bar{j} k}'$ as
\bea
\label{scalar_eqn}
\sigma_2(g')
&=& {n(n-1)\over 2}e^{2u}(1-4\alpha e^{-u}|D u|^2)+2\alpha n(n-1)fe^{-u}|D u|^2
\nonumber\\
&&
+
n(n-1)f+{n(n-1)\over 2}e^{-2u}f^2-2n\alpha\mu
\nonumber\\
&&
+2\alpha n(n-1)e^{-u}(\Delta f-2{\rm Re}(g^{j\bar k}f_ju_{\bar k})).
\eea
It follows that equations (\ref{FY1}) and (\ref{scalar_eqn}) are equivalent when $\alpha \neq 0$. Here as in the rest of the paper, we denote by $D$ the covariant derivative with respect to the given metric $g_{\bar kj}$.
Furthermore, as in \cite{FY},
we impose a normalization condition on a solution $u$. Let $\beta = {n \over n-1}$, and $\gamma = {4 \over \beta -1}$. For $A \ll 1$, we impose
\be 
\label{A_normalize}
\left( \int_X e^{- \gamma u} \right)^{1 \over \gamma} = A.
\ee
The ellipticity condition for equation (\ref{scalar_eqn}) is that the eigenvalues of $g_{\bar{j} k}'$ with respect to the metric $g_{\bar j k}$ should be in the $\Gamma_2$ cone,
\bea
\Gamma_2=\{\lambda'\in {\bf R}^n;
\ \sigma_1(\lambda')>0, \ \sigma_2(\lambda')>0\}.
\eea
Moreover, we remark that $g'\in \Gamma_2$ implies that $\tilde g_{\bar{j} k}>0$ by relation (\ref{relationg}).

\smallskip
The equation (\ref{scalar_eqn}) fits in the framework of complex Hessian equations on closed manifolds, which have been studied extensively by many authors in recent years, see for example, \cite{B2, DK, DK1, HMW, KN, Lu, LN, Sun, Gabor, GTW, Zhang1, Zhang}. However, in comparison with previous works, (\ref{scalar_eqn}) has two new difficulties. The first difficulty is the dependence on the gradient of the right hand side of the equation. This causes some trouble when attempting to obtain a $C^2$ estimate. The second difficulty is the possible degeneracy of the equation. It is easy to see that even for the ideal case $f=\mu=0$ in equation (\ref{scalar_eqn}), the right hand side might be zero. Therefore, to get smooth solutions, one needs to show that it is not degenerate under certain conditions on $A$. See \S4 and \S7 for more discussions of this particular difficulty.
Before moving to next subsection, we want to emphasize that these two difficulties occur when $\alpha>0$ in equation (\ref{scalar_eqn}). If $\alpha<0$, the behavior of the equation is quite different and Fu-Yau studied the $n=2$ case in \cite{FY1}. We will investigate the higher dimensional case in other work.

\subsection{The linearization $F^{j\bar k}$ of $\sigma_2(g')$}

We can view the Fu-Yau equation (\ref{scalar_eqn}) as a complex Hessian equation of $\sigma_2$ type, with a right hand side depending on $Du$. In accordance with standard notation in partial differential equations, we also denote $\sigma_2(g')$ by $F$, viewed as a function of $u$, $Du$, and $D\bar Du$. In particular, $F^{j\bar k}\equiv \p F/\p g'_{\bar kj}$, and the linearization of $\sigma_2(g')$ is given by 
\bea
\delta\sigma_2(g')=F^{j\bar k}\delta g_{\bar kj}'.
\eea
We shall need explicit formulas for $F^{j\bar k}$, and for the operator $2n\alpha F^{j\bar k}D_jD_{\bar k}$ acting on $u$, the gradient $Du$ of $u$, the square $|Du|^2$ of the gradient, and the complex hessian $D_pD_{\bar q}u$. 

\smallskip
We summarize briefly here our notations and conventions. The Hermitian form $\o$ defined by a K\"ahler metric $g_{\bar kj}$ is given by $\o=ig_{\bar kj}dz^j\wedge d\bar z^k$. The Chern unitary connection with respect to the metric $\o$ is denoted by $D_{\bar j}={\p\over \p \bar z^j}$,
$D_jV^p=     g^{p \bar q}\p_j(g_{\bar q m}V^m)$, and the curvature tensor is defined by
\bea
[D_{\bar k}, D_j]V^m=-R_{\bar kj}{}^m{}_pV^p.
\nonumber
\eea
The Ricci curvature $R_{\bar kj}$ is given by $R_{\bar kj}=R_{\bar kj}{}^m{}_m$. Given a second Hermitian tensor $g_{\bar km}'$, the relative endomorphism $h_j{}^k$ 
from $g_{\bar km}'$ to $g_{\bar km}$ is defined by
\bea
h_j{}^k=g^{k\bar m}g_{\bar m j}'.
\nonumber
\eea 
Writing $\sigma_2(g')=(({\rm Tr} \,h)^2-{\rm Tr}\,h^2)/2$, we readily find
\bea
F^{j\bar k}=g^{j\bar p} g^{q\bar k}\tilde g_{\bar pq}
\eea
where $\tilde g_{\bar pq}$ is the metric introduced in (\ref{C0_ellipticity}). In particular
$F^{j\bar k}g_{\bar kj}=(n-1){\rm Tr}\,h$, and hence
\bea \label{FDDu}
2n\alpha F^{j\bar k}D_jD_{\bar k}u
&=&
F^{j\bar k}g_{\bar kj}'-(e^u+fe^{-u})F^{j\bar k}g_{\bar kj}
\nonumber\\
&=&
2 F-(n-1)(e^u+fe^{-u}){\rm Tr}\,h.
\eea
Next, the variational formula for $\sigma_2(g')$ implies 
\bea
\label{dsigma}
\p_p F=F^{j\bar k}D_pg_{\bar kj}'. 
\eea
Substituting in the definition of $g_{\bar kj}'$, we obtain the following formula for
$2n\alpha F^{j\bar k}D_jD_{\bar k} (D_pu)$,
\bea
\label{3D1}
2n\alpha F^{j\bar k}D_jD_{\bar k}(D_pu)=
\p_p F
-
(n-1){\rm Tr}\,h\, \p_p(e^u+fe^{-u}).
\eea
Similarly, we find
\bea
\label{3D2}
2n\alpha F^{j\bar k}D_jD_{\bar k}(D_{\bar p}u)
=\p_{\bar p} F
-
(n-1){\rm Tr}\,h\, \p_{\bar p}(e^u+fe^{-u})
+2n\alpha
\tilde g_{\bar\ell m}R_{\bar{p}}{}^{\bar \ell m\bar q}D_{\bar q}u
\eea
where $R_{\bar{p}}{}^{\bar \ell m\bar q}$ is the curvature of metric $\o$. This additional curvature term resulted from the commutation of covariant derivatives $D_{\bar p}D_j$ and $D_jD_{\bar p}$ when acting on $D_{\bar k}u$. 
It is now easy to deduce $2n\alpha F^{j\bar k}D_jD_{\bar k}|Du|^2$. 
Introduce the notation
\bea
|DDu|_{Fg}^2=F^{j\bar k}g^{\ell\bar m}D_{j\ell}uD_{\bar k\bar m}u,
\quad
|D\bar Du|_{Fg}^2=F^{j\bar k}g^{\ell\bar m}D_{j\bar m}uD_{\ell \bar k}u.
\eea
Then
\bea
2n\alpha F^{j\bar k}D_jD_{\bar k}|Du|^2
&=&
2n\alpha g^{\ell\bar m}F^{j\bar k}(D_jD_{\bar k}D_{\ell}u \,D_{\bar m}u
+
D_{\ell}u \,D_j D_{\bar k}D_{\bar m}u)
\nonumber\\
&&
+2n\alpha (|DDu|_{Fg}^2+|D\bar D u|_{Fg}^2)
\eea
and hence, in view of the formulas (\ref{3D1}) and (\ref{3D2}),
\bea \label{DD|Du|^2}
2n\alpha F^{j\bar k}D_jD_{\bar k}|Du|^2
&=&
g^{\ell\bar m}(\p_\ell F \p_{\bar m}u+\p_{\bar m}F \p_\ell u)
+2n\alpha \tilde g_{\bar \ell m}\p_puR^{m\bar\ell p\bar q}\p_{\bar q}u
\nonumber\\
&&
-
(n-1){\rm Tr}\,h\,g^{\ell\bar m}(\p_{\bar m}(e^u+fe^{-u})\p_\ell u+\p_\ell(e^u+fe^{-u})\p_{\bar m}u)
\nonumber\\
&&
+2n\alpha (|DDu|_{Fg}^2+|D\bar D u|_{Fg}^2).
\eea
Finally, the operator $2n\alpha F^{j\bar k}D_jD_{\bar k}$ acting on the Hessian $D_pD_{\bar q}u$ can be obtained in a similar way from differentiating the equation (\ref{dsigma}) again, giving
\bea
F^{j\bar k}D_pD_{\bar q}g_{\bar kj}'
=
\p_p\p_{\bar q}F
-
D_p({\rm Tr}\,h)D_{\bar q}({\rm Tr}\,h)
+
D_j h^j{}_p D_{\bar k}h_{\bar q}{}^{\bar k}.
\eea
We can extract the term $D_pD_{\bar q}D_jD_{\bar k} u$ from the left-hand side. Permuting the order of differentiation, we find
\bea \label{DDDDu}
&&2n \alpha F^{j \bar{k}}  D_j D_{\bar{k}} D_p D_{\bar{q}} u \\&=& 2 n \alpha F^{j \bar{k}} D_p D_{\bar{q}} D_j D_{\bar{k}} u  + 2 n \alpha \left(F^{j \bar{k}} R_{\bar{q} j \bar{k}}{}^{\bar{a}} u_{\bar{a} p} - F^{j \bar{k}} R_{\bar{q} p \bar{k}}{}^{\bar{a}} u_{\bar{a} j}\right) \nonumber\\
&=& F^{j \bar{k}} D_p D_{\bar{q}} g'_{\bar{k} j} - \left(e^u -f e^{-u}\right) u_{\bar{q} p} F^{j \bar{k}} g_{\bar{k} j} - \left(e^u +f e^{-u}\right) u_p u_{\bar{q}} F^{j \bar{k}} g_{\bar{k} j} \nonumber\\
&&+ 2e^{-u} \Re\left(u_p f_{\bar{q}}\right)F^{j \bar{k}} g_{\bar{k} j} - f_{\bar{q} p} e^{-u} F^{j \bar{k}} g_{\bar{k} j} + 2 n \alpha \left(F^{j \bar{k}} R_{\bar{q} j \bar{k}}{}^{\bar{a}} u_{\bar{a} p} - F^{j \bar{k}} R_{\bar{q} p \bar{k}}{}^{\bar{a}} u_{\bar{a} j}\right) \nonumber\\
&=& \p_p\p_{\bar q} F - D_p({\rm Tr}\,h)D_{\bar q}({\rm Tr}\,h) + D_j h^j{}_p D_{\bar k}h_{\bar q}{}^{\bar k} + 2 n \alpha \left(F^{j \bar{k}} R_{\bar{q} j \bar{k}}{}^{\bar{a}} u_{\bar{a} p} - F^{j \bar{k}} R_{\bar{q} p \bar{k}}{}^{\bar{a}} u_{\bar{a} j}\right) \nonumber\\
&& + \left\{- (e^u -f e^{-u}) D_p D_{\bar{q}}u - (e^u +f e^{-u}) D_p u D_{\bar{q}} u+ 2e^{-u} \Re(u_p f_{\bar{q}}) - e^{-u} D_p D_{\bar{q}} f \right\} (n-1) {\rm Tr}\,h. \nonumber
\eea
All these formulas are quite general. For the specific Fu-Yau equation, we can substitute the right hand side of equation (\ref{scalar_eqn}) for $F=\sigma_2(g')$, as we shall do in sections \S 4 and \S 7.

\section{The $C^0$ Estimate}
\setcounter{equation}{0}

The following $C^0$ estimate holds:

\begin{theorem}
Let $(X,\o)$ be a compact K\"ahler manifold of dimension $n$ with ${\rm Vol}(X,\o)=1$. Let $u$ be a solution of (\ref{scalar_eqn}) under ellipticity condition (\ref{C0_ellipticity}) and normalization condition (\ref{A_normalize}). 
%Let $C_0$, $C_1$ be a constants which depends on $(X,\o)$, $\alpha$, $\|f\|_{C^2}$, and $\|\mu\|_{L^\infty}$. 
Then, for $A < 1$, there exists a constant $C_0$ depending only on $(X,\o)$, $\alpha$, $\|f\|_{C^2}$, and $\|\mu\|_{L^\infty}$ such that 
\be
e^{-\inf u} \leq C_0 A.
\ee
Furthermore, if $A$ is chosen small enough such that $C_0 A <1$, then there is a constant $C_1$ depending also only on $(X,\o)$, $\alpha$, $\|f\|_{C^2}$, and $\|\mu\|_{L^\infty}$ such that  
\be
e^{\sup u} \leq C_1 A^{-1}.
\ee 
\end{theorem}

{\it Proof.}
We proceed by Moser iteration. First, we define the Hermitian form corresponding to $\tilde{g}_{\bar{j} i}$:
\be
\tilde{\o} = (n-1)(e^u + f e^{-u}) \o  + 2n \alpha ( (\Delta u)\o - i \partial \bar{\partial} u) >0.
\ee
Let $k \geq 2$. The starting point is to compute the quantity
\be
\int_X i \ddb (e^{-ku}) \wedge \tilde{\o} \wedge \o^{n-2}
\ee
in two different ways. On one hand, by the definition of $\tilde{\o}$ and Stokes' theorem, we have
\be
\int_X i \ddb (e^{-ku}) \wedge \tilde{\o} \wedge \o^{n-2} = \int_X (n-1) (e^u + f e^{-u}) i \ddb (e^{-ku}) \wedge \o^{n-1} + 2 n \alpha \int_X (\D u) i \ddb (e^{-ku}) \wedge \o^{n-1}.
\ee
Using the volume form $\o^n \over n!$, we compute
\bea \label{int_identity1}
& \ &
{1 \over (n-1)!} \int_X i \ddb (e^{-ku}) \wedge \tilde{\o} \wedge \o^{n-2} \\ &=& \int_X (n-1) (e^u + f e^{-u}) \D (e^{-ku}) + 2 n \alpha \int_X (\D u) \D (e^{-ku}) \nonumber \\
&=& k^2 (n-1) \int_X (e^u + f e^{-u}) e^{-ku} | D u|^2 - k(n-1) \int_X (e^u + f e^{-u}) e^{-ku} \D u \nonumber \\
&& + 2 k^2 n \alpha \int_X e^{-ku} \D u | D u|^2 - 2k n \alpha \int_X e^{-ku} (\D u)^2\nonumber.
\eea
On the other hand, using equation (\ref{FY1}), we obtain
\bea
%&\ &
&&\int_X i \ddb (e^{-ku}) \wedge \tilde{\o} \wedge \o^{n-2} \\ &=& k^2 \int_X e^{-ku} i \p u \wedge \bar{\p} u \wedge \tilde{\o} \wedge \o^{n-2}  -k \int_X (n-1)  e^{-ku} (e^u + f e^{-u}) i \ddb u \wedge \o^{n-1} \nonumber \\
&& -2k n \alpha \int_X  e^{-ku} \D u \ i \ddb u \wedge \o^{n-1} + 2k n \alpha \int_X  e^{-ku} i \ddb u \wedge i \ddb u \wedge \o^{n-2}  \nonumber \\
&=&  k^2 \int_X e^{-ku} i \p u \wedge  \bar{\p} u \wedge \tilde{\o} \wedge \o^{n-2}  -k \int_X (n-1)  e^{-ku} (e^u + f e^{-u}) i \ddb u \wedge \o^{n-1} \nonumber \\
&& -2k n \alpha \int_X  e^{-ku} \D u \ i \ddb u \wedge \o^{n-1} - 2k(n-2)! \int_X e^{-ku} \mu {\o^n \over n!} -2k \int_X  e^{-ku} i \ddb (e^u - f e^{-u}) \wedge \o^{n-1}.\nonumber
\eea
Expanding out terms and using the definition of $\tilde{\o}$ yields
\bea \label{int_identity2}
&\ &
{1 \over (n-1)!}\int_X i \ddb (e^{-ku}) \wedge \tilde{\o} \wedge \o^{n-2} \\
&=& k^2 (n-1) \int_X e^{-ku}(e^u + f e^{-u}) | D u|^2   + k^2(2n \alpha)  \int_X e^{-ku} \D u | D u|^2 \nonumber \\
&& -{k^2 (2 n \alpha) \over (n-1)!} \int_X e^{-ku} i \p u \wedge \bar{\p} u \wedge i \ddb u \wedge \o^{n-2}  -k(n-1) \int_X  e^{-ku} (e^u + f e^{-u}) \D u \nonumber \\
&& -2k n \alpha \int_X  e^{-ku} (\D u)^2 -{2k \over n-1} \int_X e^{-ku} \mu  -2k \int_X e^{-(k-1)u} (| D u|^2 + \D u) -2k \int_X  e^{-(k+1)u} f \D u \nonumber\\
&&+ 2k \int_X e^{-(k+1)u} \D f  +2k \int_X  e^{-(k+1)u} f | D u|^2 - 4k \int_X e^{-(k+1)u} Re(g^{j \bar k}f_j u_{\bar k})  .\nonumber
\eea
We now equate (\ref{int_identity1}) and (\ref{int_identity2}) and cancel repeating terms.
\bea
0 &=& -{k n \alpha \over (n-1)!} \int_X e^{-ku} i \p u \wedge \bar{\p} u \wedge i \ddb u \wedge \o^{n-2}  - \int_X e^{-(k-1)u} | D u|^2  \\
&& - {1 \over n-1} \int_X e^{-ku} \mu-  \int_X e^{-(k-1)u} \D u   +\int_X  e^{-(k+1)u} \D f + \int_X  e^{-(k+1)u} f | D u|^2 \nonumber\\
&& - 2 \int_X e^{-(k+1)u} Re(g^{j \bar k}f_j u_{\bar k})  - \int_X  e^{-(k+1)u} f \D u .\nonumber
\eea
Integration by parts gives
\bea \label{integrals_byparts}
0 &=& -{k n \alpha \over (n-1)!} \int_X e^{-ku} i \p u \wedge \bar{\p} u \wedge i \ddb u \wedge \o^{n-2}  - {1 \over n-1}\int_X e^{-ku} \mu  - k \int_X e^{-(k-1)u} | D u|^2  \nonumber \\
&&+\int_X  e^{-(k+1)u} \D f  -k \int_X  e^{-(k+1)u} f | D u|^2 - \int_X e^{-(k+1)u} g^{j \bar k}f_j u_{\bar k}.
\eea
One more integration by parts yields the following identity:
\bea \label{integrals_simplified}
&& k \int_X e^{-ku} | D u|^2 (e^u + f e^{-u}) \\&=& -{k n \alpha \over (n-1)!} \int_X e^{-ku} i \p u \wedge \bar{\p} u \wedge i \ddb u \wedge \o^{n-2} - {1 \over n-1}\int_X e^{-ku} \mu   + (1- {1 \over k+1})\int_X  e^{-(k+1)u} \D f . \nonumber 
\eea
We now estimate the first term on the right hand side. At a point $p \in X$, choose coordinates such that $g_{\bar{k} j} = \delta_{k j}$ and $u_{\bar{k} j}$ is diagonal. From the condition $\tilde{g}>0$, we see that $\tilde{g}_{\bar{k} k} = (n-1)(e^u + f e^{-u}) + 2 n \alpha (\D u - u_{\bar{k} k}) >0$
%\be
%\tilde{g}_{\bar{k} k} = (n-1)(e^u + f e^{-u}) + 2 n \alpha (\D u - u_{\bar{k} k}) >0,
%\ee
at $p$. We compute
\bea \label{pos_estimate}
 i \p u \wedge \bar{\p} u \wedge i \ddb u \wedge \o^{n-2} &=& (n-2)! \sum_i |u_i|^2 (\D u - u_{\bar{i} i}) {\o^n \over n!}\nonumber \\
&>& - { (n-1)! \over 2 n \alpha} | D u|^2 (e^u + f e^{-u}) {\o^n \over n!} .
\eea
Using this inequality in (\ref{integrals_simplified}), we obtain
\bea
 {k \over 2} \int_X e^{-ku} | D u|^2 (e^u + f e^{-u}) &\leq& - {1 \over n-1} \int_X e^{-ku} \mu + \bigg( 1- {1 \over k+1} \bigg)\int_X  e^{-(k+1)u} \D f . \nonumber \\
\eea
Since $f \geq 0$, we can deduce the following estimate:
\be
 k \int_X e^{-(k-1)u} | D u|^2 \leq C \left( \int_X e^{-ku} + \int_X  e^{-(k+1)u} \right) . 
\ee
Therefore, for $k \geq 1$, we have
\be \label{reverse_sobolev}
 \int_X | D e^{-{k \over 2}u}|^2 \leq C k \left( \int_X e^{-(k+1)u} + \int_X  e^{-(k+2)u} \right) . 
\ee
To obtain a $C^0$ estimate, we use the method of Moser iteration as done in \cite{FY}. We set $\beta = { n \over (n-1)}$. The Sobolev inequality gives us 
\be
 \left( \int_X | e^{-{k \over 2}u}|^{2 \beta} \right)^{1 \over \beta} \leq C  \left( \int_X | e^{-{k \over 2}u}|^2  + \int_X |D e^{-{k \over 2}u}|^2 \right) . 
\ee
Combining the Sobolev inequality with (\ref{reverse_sobolev}) yields
\be
 \left( \int_X  e^{-k \beta u} \right)^{1 \over \beta} \leq C k \left( \int_X e^{-ku}+ \int_X e^{-(k+1)u} + \int_X  e^{-(k+2)u} \right) . 
\ee
Applying H\"older's inequality, we get
\be \label{integral_estimate}
 \left( \int_X  e^{-k \beta u} \right)^{1 \over \beta} \leq C k \bigg\{ \left(\int_X e^{-(k+2)u} \right)^{k \over k+2}+ \left( \int_X e^{-(k+2)u} \right)^{k \over k+1} + \int_X  e^{-(k+2)u} \bigg\} . 
\ee
For this inequality to be useful, we need to take $k$ large enough so that $k \b \geq k+2$. In order to proceed with the iteration, we consider two cases.
\newline
\newline
\textit{Case 1:} For all $k \geq \gamma = {4 \over \beta - 1}$, we have $\int_X e^{-ku} \leq 1$. In this case, for each $k \geq \gamma$, (\ref{integral_estimate}) gives us
\be 
 \left( \int_X  e^{-k \beta u} \right)^{1 \over \beta} \leq C k  \left(\int_X e^{-(k+2)u} \right)^{k \over k+2} . 
\ee
Using H\"older's inequality, we also have
\be \label{holder_split}
\int_X e^{-(k+2)u} = \int_X \left(e^{-u}\right)^{\beta \gamma \over 2} \left(e^{-u}\right)^{2k- \gamma \over 2} \leq \left( \int_X e^{- k \beta u} \right)^{\gamma \over 2 k} \left( \int_X e^{- k u} \right)^{1 -{\gamma \over 2 k}}. 
\ee
Therefore
\be 
 \left( \int_X  e^{-k \beta u} \right)^{1 \over \beta} \leq C k  \left(\int_X e^{-k \beta u} \right)^{\gamma \over 2(k+2)} \left(\int_X e^{-k u} \right)^{2k - \gamma \over 2(k+2)}. 
\ee
By regrouping and using the identity $\gamma \beta = 4 + \gamma$, we obtain

\be 
 \left( \int_X  e^{-k \beta u} \right)^{1 \over \beta} \leq \left(C k\right)^{2(k+2) \over 2k -\gamma} \int_X e^{-k u}.
\ee
Since $k \geq \gamma \geq 4$, we have ${2(k+2) \over 2k - \gamma} \leq {2(k+2) \over k} \leq 3$.
%\be 
%{2(k+2) \over 2k - \lambda} \leq {2(k+2) \over k} \leq 3.
%\ee
Thus
\be \label{case1_est}
\|e^{-u}\|_{L^{k \beta}} \leq (C k)^{3/k} \|e^{-u}\|_{L^k},
\ee
for $k \geq\gamma$. We iterate this estimate and conclude
\be 
e^{- \inf u} = \|e^{-u}\|_{L^\infty} \leq C \|e^{-u}\|_{L^{\gamma}} = C A.
\ee
\medskip
%\newline
\textit{Case 2:} There exists a $k_0 > \gamma$ such that $\int_X e^{-k_0 u} > 1$. In this case, using ${\rm Vol}(X,\o)=1$ and H\"older's inequality, we have  $\int_X e^{-k u} > 1$ for all $k \geq k_0$. After possibly increasing $k_0$, we take $k \geq k_0 \geq \gamma \beta > \gamma$. From (\ref{integral_estimate}) and (\ref{holder_split}), we have
\bea
 \left( \int_X  e^{-k \beta u} \right)^{1 \over \beta} \leq C k \int_X e^{-(k+2) u} \leq Ck \left( \int_X  e^{-k \beta u} \right)^{\gamma \over 2k} \left( \int_X  e^{-k u} \right)^{1 - {\gamma \over 2k}}.
\eea
After rearranging, we obtain
\be \label{case2_est}
\|e^{-u}\|_{L^{k \beta}} \leq (C k)^{2 \over 2k - \gamma \beta} \|e^{-u}\|_{L^k}^{2k - \gamma \over 2k - \gamma \beta}.
\ee
Since we assume $k \geq \gamma \beta$, we conclude
\be 
\|e^{-u}\|_{L^{k \beta}} \leq (C k)^{2 \over k} \|e^{-u}\|_{L^k}^{2k - \gamma \over 2k - \gamma \beta}.
\ee
We set
\be
\Theta(\nu) = {1 \over \beta} \left( {2 k_0 \beta^\nu - \gamma \over 2k_0 \beta^{\nu-1} - \gamma} \right).
\ee
For $i, j \in \N$ with $j \geq 0$ and $i \geq j$, we have
\be
\prod_{\nu = j}^i \Theta(\nu) = {2k_0 - {\gamma \over \beta^i} \over 2k_0 - {\gamma \over \beta^{j-1}}} \leq {2 k_0 \over 2k_0 -\gamma \beta} \leq 2.
\ee
 Therefore we can iterate our estimate in the following way:
\bea
\|e^{-u}\|_{L^{k_0 \beta^{i+1}}} &\leq& (C k_0 \beta^i)^{2 \over k_0 \beta^i} \|e^{-u}\|_{L^{k_0 \beta^i}}^{\Theta(i)}\leq \left( \prod_{j=0}^i (C k_0 \beta^j)^{{2 \over k_0 \beta^j} \prod_{\nu = j+1}^i \Theta(\nu)}  \right) \|e^{-u}\|_{L^{k_0}}^{\prod_{\nu=0}^i \Theta(\nu)} \nonumber\\
&\leq&  \left( \prod_{j=0}^\infty (C k_0 \beta^j)^{{4 \over k_0 \beta^j}}  \right) \|e^{-u}\|_{L^{k_0}}^2  \leq C \|e^{-u}\|_{L^{k_0}}^{2}.
\eea
As we let $i \rightarrow \infty$, we have
\be
\|e^{-u}\|_{L^\infty} \leq C \|e^{-u}\|_{L^{k_0}}^2.
\ee
We would like to estimate $\|e^{-u}\|_{L^{k_0}}$ in terms of $\|e^{-u}\|_{L^{\gamma}}$. Starting from (\ref{integral_estimate}), we can follow either case 1 or case 2, depending on the size of $\int e^{-(k+2)u}$. We then arrive at estimate (\ref{case1_est}) or (\ref{case2_est}):
\be
\|e^{-u}\|_{L^{k_0}} \leq C \|e^{-u}\|_{L^{k_0 \over \beta}}, \ {\rm or} \ \|e^{-u}\|_{L^{k_0}} \leq C \|e^{-u}\|_{L^{k_0 \over \beta}}^r, \ {\rm for \ some} \ r >1.
\ee
By repeating this process finitely many times, we can control $\|e^{-u}\|_{L^{k_0}} \leq C \|e^{-u}\|_{L^{\gamma}}^a$
%\be
%\|e^{-u}\|_{L^{k_0}} \leq C \|e^{-u}\|_{L^{\gamma}}^a,
%\ee
for some $a \geq 1$. Since $\|e^{-u}\|_{L^{\gamma}} = A <1$, we have
\be
e^{- \inf u} = \|e^{-u}\|_{L^\infty} \leq C \|e^{-u}\|_{L^{\gamma}}^{2a} \leq C A.
\ee
\newline
To control the supremum of $u$, we replace $k$ with $-k$ in (\ref{integrals_simplified}). Then, for $k \neq 1$,
%This gives the following identity for $k >1$ 
\bea
 &&k \int_X e^{ku} | D u|^2 (e^u + f e^{-u}) \\&=& -{k n \alpha \over (n-1)!} \int_X e^{ku} i \p u \wedge \bar{\p} u \wedge i \ddb u \wedge \o^{n-2}  + {1 \over n-1} \int_X e^{ku} \mu  - (1- {1 \over 1-k})\int_X  e^{(k-1)u} \D f \nonumber.  
\eea
Proceeding as before in the case of the infimum estimate, we can use (\ref{pos_estimate}) to derive the following estimate for any $k$ greater than a fixed number greater than $1$
\be
 k \int_X e^{(k+1)u} | D u|^2 \leq C \left( \int_X e^{ku} + \int_X e^{(k-1)u} \right).
\ee
Thus for $k \geq 2 \beta$, we can estimate
\be
 \int_X | D e^{{k \over 2}u}|^2 \leq C k \left( \int_X e^{(k-1)u} + \int_X e^{(k-2)u} \right).
\ee
Since $e^{- \inf u} = C A \ll 1$, we can conclude
\be \label{sup_iter}
 \int_X | D e^{{k \over 2}u}|^2 \leq C k \int_X e^{ku},
\ee
for $k \geq 2 \beta$. The Sobolev inequality yields
\be
\left( \int_X e^{k \beta u} \right)^{1/\beta} \leq C k \int_X e^{ku}.
\ee
By iterating this estimate, we have
\be \label{sup_iter}
e^{\sup u} \leq C \|e^u\|_{L^{2 \beta}}.
\ee
To complete the supremum estimate, we need another inequality. Setting $k=-1$ in (\ref{integrals_byparts}), we have 
\bea
\int_X e^{2u} | D u|^2 &=& -{ n \alpha \over (n-1)!} \int_X e^{u} i \p u \wedge  \bar{\p} u \wedge i \ddb u \wedge \o^{n-2} \nonumber \\
&& + {1 \over n-1} \int_X e^{u} \mu  - \int_X \D f - \int_X  f | D u|^2 + \int_X g^{j \bar k}f_j u_{\bar k}.
\eea
We estimate the first term on the RHS by using (\ref{pos_estimate}), and since $e^u \geq 1$ we obtain
\bea
\int_X e^{2u}| D u|^2 \leq C \left( \int_X e^u + \int_X |u| + 1 \right) \leq C \left( \int_X e^u \right).
\eea
Therefore
\be \label{L2_grad_est}
\int_X | D e^u |^2 \leq C \left( \int_X e^u \right).
\ee
Either by using this estimate, or using a scaling argument, one can obtain from (\ref{sup_iter}) that
\be \label{sup_estimate}
e^{\sup u} \leq C \|e^u\|_{L^2},
\ee
so the objective now is to control $\|e^u\|_{L^2}$. Consider the set $U = \{x : e^u \leq {2 \over A} \}$.
%\be
%U = \{x : e^u \leq {2 \over A} \}.
%\ee
We have
\bea
A^\gamma &=& \int_U e^{-\gamma u} +  \int_{X \backslash U} e^{-\gamma u} \leq e^{-\gamma \inf u} |U| + {A^\gamma \over 2^\gamma} (1 - |U|) \nonumber\leq  \left( C A^\gamma - {A^\gamma \over 2^\gamma} \right) |U| + {A^\gamma \over 2^\gamma}.
\eea
Therefore,
\be
|U| \geq \frac{1 - {1 \over 2^\gamma}}{C - {1 \over 2^\gamma}} := \delta >0.
\ee
To estimate the $L^2$ norm of $e^u$, we follow the argument from Tosatti-Weinkove \cite{TW}. Let $\psi = e^u$, and let $\underline{\psi} := \int_X \psi$. By the Poincar\'e inequality and (\ref{L2_grad_est}),
\be \label{poincare}
\| \psi - \underline{\psi}\|_{L^2} \leq \| D \psi \|_{L^2} \leq C \| \psi \|_{L^1}^{1/2}.
\ee
We compute
\bea
\delta \ \underline{\psi} \leq \int_U | \underline{\psi}| \leq \int_U |\psi - \underline{\psi}| + |\psi|\leq   \int_X  |\psi - \underline{\psi}| + {2 \over A} |U| \leq  {C \over A}( 1 + \|\psi - \underline{\psi}\|_{L^1}).
\eea
Using the previous estimate and (\ref{poincare}), it is now easy to obtain an $L^1$ estimate:
\bea
\| \psi\|_{L^1} &\leq&  \|\psi - \underline{\psi}\|_{L^1} + \|\underline{\psi}\|_{L^1}\leq   CA^{-1}( 1 + \|\psi - \underline{\psi}\|_{L^1}) \nonumber\\
&\leq&  CA^{-1}( 1 + \|\psi - \underline{\psi}\|_{L^2})\leq  CA^{-1}( 1 + \| \psi \|_{L^1}^{1/2}).
\eea
Therefore $\| \psi \|_{L^1}$ is under control, and by (\ref{poincare}), we can deduce that $\|\psi\|_{L^2} \leq CA^{-1}$. By (\ref{sup_estimate}), we have
\be
e^{\sup u} \leq CA^{-1}.
\ee

\section{The $C^1$ Estimate}
\setcounter{equation}{0}

\par
As we mentioned previously, the {\it a priori} gradient estimate is easy due to the special structure of the right hand side of equation (\ref{scalar_eqn}). %The $C^1$ estimate follows easily from the equation. 
Define the constant $\kappa_c$ by
\be\label{definekc}
\kappa_c = {n (n-1) \over 2}.
\ee
Our equation (\ref{scalar_eqn}) is
\bea \label{kappa_equation}
e^{-2u} F
&=&
\kappa_c - 4 \alpha \kappa_c \bigg \{  e^{-u} |D u|^2 - f e^{-3u} |D u|^2  + e^{-3u} (g^{i \bar{k}}f_i u_{\bar{k}} + g^{i \bar{k}} f_{\bar{k}} u_i) \bigg \}
\nonumber\\
&&
+ \kappa_c e^{-2u}  \bigg\{ 2f + f^2 e^{-2u} + 4 \alpha e^{-u} \Delta f \bigg\}
-2n\alpha e^{-2u} \mu.
\eea
We first estimate
\be \label{first_estimate}
0 < e^{-2u} F \leq \kappa_c - 4 \alpha \kappa_c e^{-u} |D u|^2 \bigg \{ 1 - (\|f\|_\infty +1) e^{-2u} \bigg \} + O(e^{-2u}).
\ee
For a choice of $A$ small enough, we can make $e^{-u} \leq C A \ll 1$.
%\be
%e^{-u} \leq C A \ll 1.
%\ee
It follows that
\be
% \label{gradient_bound}
e^{-u} |D u|^2 \leq C.
\ee
\begin{theorem}\label{C1estimate}
Let $u$ be a solution of (\ref{scalar_eqn}) under ellipticity condition (\ref{C0_ellipticity}) and normalization condition (\ref{A_normalize}). If $A$ is small enough, then there exists a positive constant $C$ depending on $(X,g)$, $\alpha$, $\|f\|_{C^2}$, and $\|\mu\|_{L^\infty}$ such that
\be \label{gradient_bound}
e^{-u} |D u|^2 \leq C.
\ee
\end{theorem}

\par
We observe that the present situation is different from the situation
for the standard complex Hessian equation
\begin{eqnarray}\label{chessian}
\sigma_k(g_{\bar ji} + u_{\bar ji} ) = f(z) 
\end{eqnarray}
on a compact K\"ahler manifold $(X, g)$ with $0<f(x) \in C^{\infty}(X)$, see \cite{DK, HMW, Zhang}. For the standard equation (\ref{chessian}) with non-degenerate right hand side $f(z)$, one needs to work very hard to get the gradient estimate since the upper bound of the $C^2$ estimate depends on the $C^1$ estimate. Once the gradient estimate is obtained, the non-degeneracy of $f(z)$ together with the $C^2$ upper bound imply the uniform ellipticity of the equation. In our current situation, the structure of equation (\ref{FY}) is better in the sense that it automatically gives a $C^1$ upper bound. However, in this case, the $C^1$ estimate is not good enough to give uniform ellipticity. For that purpose, we need to get a uniform positive lower bound for $e^{-2u}F$, which turns out to be equivalent to a sharper $C^1$ upper bound. From this viewpoint, the desired gradient estimate here is much more involved than in the standard case. We will continue to discuss this in \S7.

%However, to show the $C^2$ estimate, we need to get an improved $C^1$ estimate, which is equivalent to the lower bound of $\sigma_2(g')$. This is not the same situation as the standard complex Hessian equation (see \cite{DK, HMW}). 

\section{The $C^2$ Estimate}
\setcounter{equation}{0}

In this section, we derive the {\it a priori} $C^2$ estimate of equation (\ref{scalar_eqn}) under the assumption of a {\it sharp gradient estimate}. As previously mentioned, the presence of the gradient of $u$ on the right hand side brings substantial difficulties. For real Hessian equations, this problem was recently addressed by Guan-Ren-Wang \cite{GRW} under some assumptions. Here, we adapt some of their ideas to the complex setting. However, there are still some troublesome terms such as $|DDu|^2$ which cannot be handled as in the real case. This is the reason for the {\it sharp gradient estimate} assumption in our estimate. Our theorem is the following.
\begin{theorem} \label{C2estimate}
Let $u$ be a solution of (\ref{scalar_eqn}) under ellipticity condition (\ref{C0_ellipticity}) and normalization condition (\ref{A_normalize}). Suppose that for every $0< \delta < 1$, there exists an $0<A_\delta \ll 1$ such that for all $0<A \leq A_\delta$, the following bound holds:
\be \label{C2_assumption}
e^{-u} |D u|^2 \leq \delta.
\ee
Then there exists $0<A_0 \ll 1$ such that for all $0<A \leq A_0$, there holds
\be\label{C2estimate1}
{1 \over C} g \leq \tilde{g} \leq C g,
\ee
where $C$ is a constant depending on $\| u \|_{L^\infty}$, $\|D u\|_{L^\infty}$, $(X,\omega)$, $\|f\|_{C^2}$, $\|\mu\|_{C^2}$, $\alpha$, $A$.
\end{theorem}

\medskip
\par
Let $B_0$, $B_1$ be constants depending on $(X,\omega)$, $\|f\|_{C^2}$, $\|\mu\|_{C^2}$, $\alpha$. Recall that we have the following $C^0$ estimates
\be \label{C2_remindC0_estimate}
e^{-u} \leq B_0 A \ll 1, \ \ e^u \leq B_1 A^{-1}.
\ee
The estimate in the assumption (\ref{C2_assumption}) was obtained by Fu-Yau in \cite{FY} when $X$ has dimension $n=2$. The Fu-Yau estimate is rederived in  \S \ref{lower_b_section} and can be found in (\ref{sharp_grad}). Whether a Fu-Yau type gradient estimate holds for dimension $n>2$ is still unknown. From equation (\ref{kappa_equation}), one can see that such an estimate implies a lower bound for $e^{-2u}F$. For the purpose of the $C^2$ estimate, we shall take
\be \label{F_lower}
e^{-2u} F \geq {1 \over 2}.
\ee
To prove the theorem, it suffices to obtain an upper bound on the maximal eigenvalue of $g'$. The upper and lower bounds of $\tilde{g}$ will then follow from the relations between $g'$ and $\tilde{g}$ as discussed in \S $2$.
\newline
\newline
Before proceeding with the $C^2$ estimate, we state a lemma due to Guan-Ren-Wang \cite{GRW}.
\begin{lemma} \label{guan-ren-wang}
Suppose $v^{i}{}_j$ is an endomorphism such that $v\in \Gamma_2$. Then for any tensor $A_{\bar{i} j}$,
\be
- \sum_{i \neq j} A_{\bar{i} i } A_{\bar{j} j} \geq -{\left| \sigma_2^{i \bar{i}}(v) A_{\bar{i} i} \right|^2 \over |\sigma_2(v)|}.
\ee
In particular,
\be
- \sum_{i \neq j} D_k g'_{\bar{i} i} D_{\bar{k}} g'_{\bar{j} j} \geq - {|D_k F|^2 \over |F|}.
\ee
\end{lemma}
 
{\it Proof.} We reproduce the proof of the Guan-Ren-Wang inequality for completeness. Let $H = {\sigma_2 (v) \over \sigma_1 (v)}.$
%\be
%H = {\sigma_2 (v) \over \sigma_1 (v)}.
%\ee
Differentiate $\log H$ with respect to the $(p,p)$ entry to obtain
\be
{H^{p \bar{p}} \over H} = {\sigma_2^{p \bar{p}} \over \sigma_2} - {\sigma_1^{p \bar{p}} \over \sigma_1}.
\ee
Differentiate again
\be
{H^{p \bar{p}, q \bar{q}} \over H} = {H^{p \bar{p}}H^{q \bar{q}} \over H^2} + {\sigma_2^{p \bar{p}, q \bar{q}} \over \sigma_2} - {\sigma_2^{p \bar{p}}\sigma_2^{q \bar{q}} \over (\sigma_2)^2} + {\sigma_1^{p \bar{p}}\sigma_1^{q \bar{q}} \over (\sigma_1)^2} .
\ee
Since $H$ is concave, we have
\bea
0 &\geq& {|H^{p \bar{p}} A_{\bar{p} p}|^2 \over H^2} + {\sigma_2^{p \bar{p}, q \bar{q}} A_{\bar{p} p} A_{\bar{q} q} \over \sigma_2} - {|\sigma_2^{p \bar{p}} A_{\bar{p} p}|^2 \over (\sigma_2)^2} + {|\sigma_1^{p \bar{p}} A_{\bar{p} p}|^2 \over (\sigma_1)^2} \nonumber\\
&\geq&  {\sigma_2^{p \bar{p}, q \bar{q}} A_{\bar{p} p} A_{\bar{q} q} \over \sigma_2} - \left| \sigma_2^{p \bar{p}} A_{\bar{p} p} \over \sigma_2 \right|^2.
\eea
This completes the proof of the Guan-Ren-Wang inequality.
\bigskip
\par
We now proceed to the proof of the $C^2$ estimate. We shall apply the maximum principle to a function similar to the one used by Hou-Ma-Wu in \cite{HMW}. Let $M>0$ be a large constant to be determined later. Let $\sup_X |u| \leq L$. Define
\be
\psi(t) = {M \over 2n \alpha} \log \left( 1 + {t \over L} \right),
\ee
It follows that
\be \label{psi_identities}
{M \over 2n \alpha L} > \psi' > {M \over 4n \alpha L}, \ \ 2n \alpha \psi'' = - {|2n \alpha \psi'|^2 \over M}.
\ee
For small $\delta>0$ to be chosen later, we define
\be
\phi(t) = - \log (M_2 - t), \ M_2 = {17 \delta B_1 \over A}.
\ee
Note that $\phi( |D u|^2 )$ is well-defined by the assumption on gradient estimate (\ref{C2_assumption}). Indeed, we may choose $A_0 \ll 1$ depending on $\delta$ such that, for any $0<A \leq A_0$, 
\be \label{grad_upper}
|D u|^2 \leq \delta e^u \leq { \delta B_1 \over A},
\ee
and hence
\be \label{phi_upper}
\phi'( |D u|^2) \leq {A \over 16 \delta B_1}.
\ee
Furthermore, we have the lower bound
\be \label{phi_lower}
\phi'( |D u|^2) \geq {A \over 17 \delta B_1} \geq {e^{-u} \over 17 \delta B_0 B_1},
\ee
and the relationship
\be \label{phi_ode}
\phi'' = (\phi')^2 .
\ee
First, consider
\be
G_0(z,\xi)= \log ( g'_{\bar{j} k} \xi^k \bar{\xi}^j) - 2n \alpha \psi(u) + \phi( |D u|^2 ),
\ee
for $z \in X$ and $\xi \in T_z^{1,0}(X)$ a unit vector. $G_0$ is not defined everywhere, but we may restrict to the compact set where $g'_{\bar{j} k} \xi^k \bar{\xi}^j \geq 0$ and obtain an upper semicontinuous function. Let $(p, \xi_0)$ be the maximum of $G_0$. Choose coordinates centered at $p$ such that $g_{\bar{j} k} = \delta_{jk}$ and $g'_{\bar{j} k}$ is diagonal. Suppose $g'_{\bar{1} 1}$ is the largest eigenvalue of $g'$. Then $\xi_0(p) = \partial_1$, and we extend this to a local unit vector field $\xi_0 = g_{\bar{1} 1}^{-1/2} {\partial \over \partial z^1}$.
%\be
%\xi_0 = g_{\bar{1} 1}^{-1/2} {\partial \over \partial z^1}.
%\ee
Define the local function
\be
G (z)= \log (g_{\bar{1} 1}^{-1} g'_{\bar{1} 1}) - 2n \alpha \psi(u) + \phi( |D u|^2 ).
\ee
This function $G$ also attains a maximum at $p \in X$. We will compute at the point $p$. We shall be assuming that $g'_{\bar{1} 1}(p) \gg 1$, otherwise we would already have an upper bound on the maximal eigenvalue of $g'$ and the $C^2$ estimate would be complete.
\newline
\newline
Covariantly differentiating $G$ gives
\be \label{testfunction1}
G_{\bar{j}} = {(e^u + f e^{-u})_{\bar{j}} + 2n \alpha D_{\bar{j}} D_1 D_{\bar{1}} u \over {g'_{\bar{1} 1}}} + \phi' | D u|^2_{\bar{j}} - 2n \alpha \psi' u_{\bar{j}}.
\ee
Differentiating $G$ a second time and contracting with $F^{i \bar{j}}$ yields
\bea \label{testfunction2}
F^{i \bar{j}} G_{\bar{j} i} &=& {2n \alpha \over g'_{\bar{1} 1}} F^{i \bar{j}} D_i D_{\bar{j}} D_1 D_{\bar{1}} u + {(e^u - f e^{-u}) \over g'_{\bar{1} 1}} F^{i \bar{j}} u_{\bar{j} i} + {(e^u + f e^{-u}) \over g'_{\bar{1} 1}} |Du|^2_F \nonumber\\
&& -{2 e^{-u} \over g'_{\bar{1} 1}}  \Re (F^{i \bar{j}} u_i f_{\bar{j}}) + { e^{-u} \over g'_{\bar{1} 1}} F^{i \bar{j}}f_{\bar{j} i} - {|D g'_{\bar{1} 1}|^2_F \over (g'_{\bar{1} 1})^2} + \phi' F^{i \bar{j}} |D u|^2_{\bar{j} i} + \phi'' | D |D u|^2|^2_F  \nonumber\\
&& -2n \alpha \psi' F^{i \bar{j}} u_{\bar{j} i} -2n \alpha \psi'' |D u|^2_F.
\eea
Here we introduced the notation
\be
| D\chi |^2_F = F^{j \bar{k}} D_j \chi D_{\bar{k}} \chi.
\ee
We will get an estimate for $D_iD_{\bar{j}} D_1 D_{\bar{1}}u$ using our formula (\ref{DDDDu}). First, notice
\be \label{g'_vs_Tr}
g'_{\bar{1} 1} \leq {\rm Tr}\,h \leq n g'_{\bar{1} 1}.
\ee
Furthermore, since $g' \in \Gamma_2$, we can estimate for each $k$,
\be
2n \alpha |u_{\bar{k} k}| \leq |g'_{\bar{k} k
}| + |e^u + fe^{-u}| \leq C( g'_{\bar{1} 1} + 1).
\ee
Using these inequalities, we may estimate (\ref{DDDDu}) in the following way
\be
2n \alpha F^{i \bar{j}} D_i D_{\bar{j}} D_1 D_{\bar{1}} u  \geq D_1 D_{\bar{1}} F - D_1({\rm Tr}\,h)D_{\bar 1}({\rm Tr}\,h) + |D_1 g'|^2 - C (g'_{\bar{1} 1})^2 -C.
\ee
We now substitute this inequality into (\ref{testfunction2}) to obtain
\bea 
F^{i \bar{j}} G_{\bar{j} i} &\geq& {1 \over g'_{\bar{1} 1}} \left( |D_1 g'|^2 - |D_1 {\rm Tr}\,h|^2
%D_1 {\rm Tr}\,h  D_{\bar{1}} {\rm Tr}\,h  
+ D_1 D_{\bar{1}} F \right) + {(e^u - f e^{-u}) \over g'_{\bar{1} 1}} F^{i \bar{j}} u_{\bar{j} i}  - {|D g'_{\bar{1} 1} |^2_F \over (g'_{\bar{1} 1})^2}   + \phi' F^{i \bar{j}} |D u|^2_{\bar{j} i}\nonumber\\
&& + \phi'' | D |D u|^2 |^2_F -2n \alpha \psi' F^{i \bar{j}} u_{\bar{j} i}  -2n \alpha \psi'' | D u|^2_F -C g'_{\bar{1} 1} - C.
\eea
We have the identity
\bea
2n \alpha F^{i \bar{j}} u_{\bar{j} i} &=& F^{i \bar{j}} g'_{\bar{j} i} -(e^u +fe^{-u}) F^{i \bar{j}} g_{\bar{j} i} = 2F -(e^u +fe^{-u}) (n-1) {\rm Tr}\,h.
\eea
Note that by estimates (\ref{C2_remindC0_estimate}), we have that $(n-1)e^u \geq 1$ for small enough choice of $A$. Using this fact with (\ref{g'_vs_Tr}), we obtain
\bea 
0 &\geq& {1 \over g'_{\bar{1} 1}} \left( |D_1 g'|^2 - D_1 {\rm Tr}\,h  D_{\bar{1}} {\rm Tr}\,h  + D_1 D_{\bar{1}} F \right)  - {|D g'_{\bar{1} 1} |^2_F \over (g'_{\bar{1} 1})^2}  + \phi' F^{i \bar{j}} |D u|^2_{\bar{j} i}  \nonumber\\
&&+ \phi'' | D |D u|^2 |^2_F + (\psi'-C) g'_{\bar{1} 1} -2n \alpha \psi'' | D u|^2_F  - C.
\eea
We now compute the term involving $\phi'$. By (\ref{DD|Du|^2}), we have
\bea \label{A2}
2n \alpha F^{i \bar{j}} |D u|^2_{\bar{j} i}
&\geq& -2 |D u| |D F| - C g'_{\bar{1} 1} + 2n \alpha( |D \bar{D} u|^2_{Fg} + |D D u|^2_{Fg}) - C.
\eea
Therefore
\bea 
0 &\geq& {1 \over g'_{\bar{1} 1}} \bigg\{ |D_1 g'|^2 - D_1 {\rm Tr}\,h  D_{\bar{1}} {\rm Tr}\,h  + D_1 D_{\bar{1}} F \bigg\} - {|D g'_{\bar{1} 1} |^2_F \over (g'_{\bar{1} 1})^2}  + \phi'|D \bar{D} u|^2_{Fg} + \phi'|D D u|^2_{Fg}  \nonumber\\
&&+ \phi'' | D |D u|^2 |^2_F - {\phi' \over n \alpha} |D u| |D F| + (\psi' - C \phi'-C) g'_{\bar{1} 1} -2n \alpha \psi'' | D u|^2_F  - C.
\eea 
Define
\be
\tau = {1 \over 1 + M}.
\ee
Using (\ref{phi_ode}), $D G(p) =0$, and (\ref{psi_identities}),
\bea
\phi'' | D |D u|^2 |^2_F &=& \sum_i F^{i \bar{i}} |\phi'|D u|^2_i|^2=\sum_i F^{i \bar{i}} \left| -{D_i g'_{\bar{1} 1} \over g'_{\bar{1} 1}}  + 2n \alpha \psi' u_i \right|^2 \nonumber\\
&\geq& \tau \sum_i F^{i \bar{i}} \left| {D_i g'_{\bar{1} 1} \over g'_{\bar{1} 1}} \right|^2 - {\tau \over 1 - \tau} \sum_i F^{i \bar{i}} \left|2n \alpha \psi' u_i \right|^2 \nonumber\\
&=& \tau \sum_i F^{i \bar{i}} \left| {D_i g'_{\bar{1} 1} \over g'_{\bar{1} 1}} \right|^2 + {\tau M \over 1 - \tau} 2n \alpha \psi'' \sum_i F^{i \bar{i}} | u_i|^2\nonumber\\
&=& \tau \sum_i F^{i \bar{i}} \left| {D_i g'_{\bar{1} 1} \over g'_{\bar{1} 1}} \right|^2 + 2n \alpha \psi'' \sum_i F^{i \bar{i}} | u_i|^2.
\eea
Thus
\bea
0 &\geq& {\phi' \over 2} |D \bar{D} u|^2_{Fg} + {\phi' \over 2}|D D u|^2_{Fg} + {1 \over g'_{\bar{1} 1}} D_1 D_{\bar{1}} F -  {\phi' \over n \alpha} |D F| |D u| \nonumber\\
&& +{1 \over g'_{\bar{1} 1}} \bigg\{ |D_1 g'|^2 - D_1 {\rm Tr}\,h  D_{\bar{1}} {\rm Tr}\,h \bigg\} -(1-\tau) {|D g'_{\bar{1} 1} |^2_F \over (g'_{\bar{1} 1})^2} \nonumber\\
&&+ {\phi' \over 2} |D \bar{D} u|^2_{Fg} + {\phi' \over 2}|D D u|^2_{Fg} +\bigg\{ {M \over 4n\alpha L}-C\phi' -C \bigg\} g'_{\bar{1} 1} - C.
\eea
Computing in coordinates and applying the Guan-Ren-Wang inequality (Lemma \ref{guan-ren-wang}) yields
\begin{eqnarray*}
 \left(|D_1 g'|^2  - D_1 {\rm Tr}\,h  D_{\bar{1}} {\rm Tr}\,h \right) &=& \sum_{i \neq j} |D_1 g'_{\bar{i} j}|^2 - \sum_{i \neq j} D_1 g'_{\bar{i} i} D_{\bar{1}} g'_{\bar{j} j}\geq \sum_{i \neq j} |D_1 g'_{\bar{i} j}|^2 -{|D_1 F|^2 \over F}.
\end{eqnarray*}
Using the definition of $g'$, we obtain
\bea
\sum_{i \neq j} |D_1 g'_{\bar{i} j}|^2 &\geq& \sum_{j > 1} |D_1 g'_{\bar{1} j}|^2=  \sum_{j > 1} |D_j g'_{\bar{1} 1} -(e^u +f e^{-u})_{j}|^2 \nonumber\\
&\geq& \left(1 - {\tau \over 2} \right)  \sum_{j > 1} |D_j g'_{\bar{1} 1}|^2 - C_M
\eea
where the last constant $C_M$ depends on $\tau$, and hence on $M$.
We therefore arrive at
\bea \label{G_ineq}
0 &\geq& {\phi' \over 2} |D \bar{D} u|^2_{Fg} + {\phi' \over 2}|D D u|^2_{Fg} + {1 \over g'_{\bar{1} 1}} D_1 D_{\bar{1}} F - {\phi' \over n \alpha} |D F| |D u| \nonumber\\
&&+ \left( 1 - {\tau \over 2} \right){1 \over g'_{\bar{1} 1}} \sum_{j > 1} |D_j g'_{\bar{1} 1}|^2  -(1-\tau) {|D g'_{\bar{1} 1} |^2_F \over (g'_{\bar{1} 1})^2}   + {\phi' \over 2} |D \bar{D} u|^2_{Fg} + {\phi' \over 2}|D D u|^2_{Fg} \nonumber\\
&& -{|D_1 F|^2 \over g'_{\bar{1} 1} F} +\bigg\{ {M \over 4n\alpha L}-C\phi' -C \bigg\} g'_{\bar{1} 1} - C_M.
\eea

At this point, it will be important to distinguish constants which depend on $A$ from those that do not. Let $B$ denote a constant depending on $(X,g)$, $\| f \|_{C^2}$, $\| \mu \|_{C^2}$, $\alpha$. As before, we use $C$ to denote a constant depending on $(X,g)$, $\| u \|_\infty$, $\| D u \|_\infty$, $\| f \|_{C^2}$, $\| \mu \|_{C^2}$, $\alpha$ and use $C_M$ to denote the constants which may also depend on $M$. We now state two lemmas.
\begin{lemma} \label{c2lemma1}
Under the non-degeneracy assumption (\ref{F_lower}) and $C^0$ estimate (\ref{C2_remindC0_estimate}), there holds
\be
\label{gradF}
|D F| \leq  B (e^{u} |D u|+e^{-u}) \left( |D D u| + |D \bar{D} u| \right) + C. 
\ee
\be \label{gradF_overF}
{|D F|^2 \over |F|} \leq B e^u \left( |D D u|^2 + |D \bar{D} u|^2 \right) +C.
\ee
\be
\label{D1D1}
|D_1 D_{\bar{1}} F| \leq B e^u  \left( |D D u|^2 + |D \bar{D} u|^2 + e^u|D u_{\bar{1} 1}| \right) +C.
\ee
\end{lemma}
\begin{lemma} \label{c2lemma2}
Let $p \in X$ be a point where $G$ attains a maximum. Assuming $g'_{\bar{1} 1}(p) \gg 1$ is large enough, then at $p$ we have
\be
{ 1 - {\tau \over 2} \over g'_{\bar{1} 1}} \sum_{j > 1} |D_j g'_{\bar{1} 1}|^2  -(1-\tau) {|D g'_{\bar{1} 1} |^2_F \over (g'_{\bar{1} 1})^2} + {\phi' \over 2} |D \bar{D} u|^2_{Fg} + {\phi' \over 2}|D D u|^2_{Fg} \geq 0.
\ee
\end{lemma}
\
\newline
Assuming these lemmas, we shall now prove the $C^2$ estimate. We may assume $g'_{\bar{1} 1} \gg 1$ is large at the point $p \in X$, otherwise we already have the desired estimate. Applying both lemmas to (\ref{G_ineq}), we have
\bea \label{c2_est_afterlemmas}
0 &\geq& {\phi' \over 2} |D \bar{D} u|^2_{Fg} + {\phi' \over 2}|D D u|^2_{Fg} - {B e^u \over g'_{\bar{1} 1}} \bigg\{ |D D u|^2 + |D \bar{D} u|^2 + e^u |D u_{\bar{1} 1}| \bigg\} -{\phi'\over n\alpha}|DF|\,|Du|
\nonumber\\
&& +\bigg\{ {M \over 4n\alpha L}-C(1+\phi') \bigg\} g'_{\bar{1} 1} - C_M.
\eea
Using $D G(p) =0$ (\ref{testfunction1}), we may estimate
\bea \label{c2_est_3rd_order}
|D u_{\bar{1} 1} | &\leq& C + |D u| g'_{\bar{1} 1} \bigg\{ \psi' + { \phi' \over 2n \alpha} (|D D u| + |D \bar{D}u|) \bigg\} \nonumber\\
&\leq &C +\left(e^{-\frac{u}{2}}  (|D D u| + |D \bar{D}u|) \right) \left( e^{\frac{u}{2}}|Du|g'_{\bar{1} 1}{ \phi' \over 2n \alpha} \right) + |D u| \psi' g'_{\bar{1} 1}\\\nonumber
&\leq& C + e^{-u} \left( |D D u|^2 + |D \bar{D} u|^2 \right) + {(\phi')^2 \over (2n\alpha)^2} {e^u |Du|^2 \over 2} (g'_{\bar{1} 1})^2 + |D u| \psi' g'_{\bar{1} 1}.
\eea
%Choose coordinates such that $g_{\bar{k} j}=\delta_{kj}$, and $u_{\bar{k} j}$, $F^{i \bar{j}}$, $g'$ are diagonal. In these coordinates,
%\be
%F^{i \bar{i}} = \sum_{k \neq i} g'_{\bar{k} k}.
%\ee
Using the estimate (\ref{c2_est_3rd_order}), the estimate (\ref{psi_identities}) for $\psi'$,
and $F^{i \bar{i}} \geq F^{1\bar{1}}$, (\ref{c2_est_afterlemmas}) becomes
\bea \label{c2_final_estimate}
0 &\geq& {1 \over g'_{\bar{1} 1}} \bigg\{ {\phi' \over 2} g'_{\bar{1} 1} F^{1 \bar{1}} - B e^u \bigg\} \left( |D D u|^2 + |D \bar{D} u|^2 \right)-{\phi'\over n\alpha}|DF|\,|Du|\nonumber\\
&&  +\bigg\{ {M \over 4n\alpha L}-C(1+ \phi' + (\phi')^2) \bigg\} g'_{\bar{1} 1} - C_M.
\eea
We shall show that for small enough $A$, we can ensure
\be
{\phi' \over 2} g'_{\bar{1} 1} F^{1 \bar{1}} - B e^u \geq 1.
\ee
Indeed, this follows from the basic fact that $g'_{\bar 1 1}F^{1\bar 1} \geq c(n) F$. Note that $g'_{\bar 11}$ is the largest eigenvalue and hence $g'_{\bar 11}\geq \frac{1}{n-1}\sigma_1(\lambda' | 1)$. We use the notation $\sigma_k(\lambda'|j)$ for the $k$-th symmetric function of $(\lambda' | j) = (\lambda_1', \cdots, \widehat{\lambda_{j}'}, \cdots, \lambda_n')\in \R^{n-1}$. For example, $\sigma_1(\lambda' | 1) = \sum_{i\neq 1} \lambda'_i= F^{1\bar 1}$. This implies
\begin{eqnarray*}
\sigma_2(\lambda')= \sigma_1(\lambda' | 1) \lambda_1 + \sigma_2(\lambda' |1) &\leq& \sigma_1(\lambda' | 1) \lambda'_1+ \frac{n-2}{2(n-1)} \sigma_1^2(\lambda' |1) \\
&\leq&  \sigma_1(\lambda' | 1) \lambda'_1+ \frac{n-2}{2} \sigma_1(\lambda' | 1) \lambda'_1,
\end{eqnarray*}
which gives the desired estimate $g'_{\bar 11}F^{1\bar 1} \geq \frac{2}{n} F$. Therefore, using (\ref{F_lower}) and (\ref{phi_lower})
\bea
g'_{\bar{1} 1}{\phi' \over 2} F^{1 \bar{1}}  - B e^u &\geq& {\phi' e^{2u} \over n} (Fe^{-2u}) - B e^u \geq e^u \left( {1 \over 2 n (17 \delta B_0 B_1)} - B \right)\geq 1,
\eea
when our parameter $A_0$ is chosen such that $\delta$ is sufficiently small and $e^u$ sufficiently large. The estimate is possible because the $B_0$, $B_1$, $B$ are independent of $A$. Now that our normalization $A \leq A_0$ has been chosen, we recall the bounds (\ref{phi_upper}) and (\ref{phi_lower}) for $\phi'$, and set $M_1={M\over 4n\alpha L}-C(1+{A\over 16\delta B_1}+({A\over 16\delta B_1})^2)$, which is positive for $M$ large enough. The inequality (\ref{c2_final_estimate}) implies
\bea
0&\geq& {1\over g_{\bar 11}'}(|DDu|^2+|D\bar Du|^2)+M_1g_{\bar 11}'-{\phi'\over n\alpha}|DF|\,|Du|-C_M
\nonumber\\
&\geq& M_1^{1\over 2}(|DDu|^2+|D\bar Du|^2)^{1\over 2}+{3\over 4}M_1g_{\bar 11}'-{\phi'\over n\alpha}|DF|\,|Du|-C_M
\nonumber\\
&\geq & {3\over 4}M_1g_{\bar 11}'-C_M
\nonumber
\eea
In the last inequality, we made use of the estimate (\ref{gradF}) for $|DF|$, and chose $M$ large enough. Thus we have established
\be
g'_{\bar{1} 1} \leq C.
\ee
This completes the proof of Theorem \ref{C2estimate}.
\newline
\newline
{\it Proof of Lemma \ref{c2lemma1}.}
Using the definition (\ref{scalar_eqn}) of $F=\sigma_2(g')$, we shall estimate $DF$ and $D_1 D_{\bar{1}} F$. The expression (\ref{scalar_eqn}) shows that $F$ is a linear combination of the expressions
\bea
e^{au}, \qquad e^{\pm u}|Du|^2,
\qquad e^{-u}Du, \qquad e^{-u}\bar Du,
\nonumber
\eea
with coefficients given by smooth functions whose derivatives of any fixed order can be bounded by constants $C$. The constant $a$ can take the values $0,\pm 1, 2$. Thus $DF$ can be bounded by a linear combination of the above expressions and their derivatives. The expressions can themselves be bounded by constants $C$, while $D(e^{au})$ can be bounded by constants $C$, and
\bea
&&
D(e^{\pm u}|Du|^2)=e^{\pm u}(DDu\cdot\bar Du+Du\cdot D\bar Du)\pm e^{\pm u}(Du)|Du|^2,
\nonumber\\
&&
D(e^{-u}Du)=e^{-u}DDu-e^{-u}Du\,Du,
\qquad
D(e^{-u}\bar Du)=e^{-u}D\bar D u-e^{-u}Du\,\bar Du.
\nonumber
\eea
All the last terms on the right hand side of each of the above three equations 
can be bounded by $C$. The estimate (\ref{gradF}) for $|DF|$ follows.

Next, we turn to $D_1D_{\bar 1}F$. For this, we view $D_{\bar 1}F$ as a linear combination
of the above expressions and their derivatives, and apply $D_1$.
In this process, we can ignore all terms bounded by
expressions of the form
\bea
e^{pu}|Du|^q(|DDu|+|D\bar D u|)+e^{ru}+|Du|^s
\nonumber
\eea
for some $p, q, r,s\geq 0$ since they can all be absorbed into
$e^u(|DDu|^2+|D\bar Du|^2)+C$ (recall that $e^u>1$ by the assumption
(\ref{C2_assumption})). Examples of such terms are the bounds obtained in (\ref{gradF})
for $|DF|$. Thus,
when the derivative $D_1$ lands on the coefficients of the linear combination giving $D_{\bar 1}F$, 
we obtain only
expressions that can be bounded by the right hand side of (\ref{gradF}) and can be ignored. This means that, to establish the bound (\ref{D1D1}),
it suffices to consider the expressions $D_1 D_{\bar 1}(e^{au})$, $D_1 D_{\bar 1}(e^{\pm u}|Du|^2)$, $D_1D_{\bar 1}(e^{-u}Du)$, and 
$D_1D_{\bar 1}(e^{-u}\bar Du)$. Modulo $O(e^{pu}|Du|^q(|DDu|+|D\bar D u|)+e^{ru}+|Du|^s)$, we can write
$D_1D_{\bar 1}(e^{au})=0$, and
\bea
D_1D_{\bar 1}(e^{-u}\bar Du)
=
e^{-u}D_1D_{\bar 1}\bar D u,
\qquad
D_1D_{\bar 1}(e^{-u}Du)=e^{-u}D_1D_{\bar 1}Du
\nonumber
\eea
which can clearly be bounded by the right hand side of (\ref{D1D1}). Similarly,
\bea
D_1D_{\bar 1}(e^{\pm u}|Du|^2)
&=&
e^{\pm u}D_1(D_{\bar 1}Du\cdot \bar Du+Du\cdot D_{\bar 1}\bar Du)
\nonumber\\
&=&
e^{\pm u}(D_1D_{\bar 1}Du\cdot \bar Du+Du\cdot D_1D_{\bar 1}\bar Du)
\nonumber\\
&&
%\qquad\qquad
+
e^{\pm u}(D_{\bar 1}Du\cdot D_1\bar Du+D_1Du\cdot D_{\bar 1}\bar Du).
\nonumber
\eea
It follows that
\bea
|D_1D_{\bar 1}(e^{\pm u}|Du|^2)|
\leq e^u(|D\bar Du|^2+|DDu|^2)+e^u|Du|\,|Du_{1\bar 1}|.
\eea
We note that by $|D u|^2 \leq e^u$ (\ref{C2_assumption}), we may estimate $|Du|$ by $e^u$ since $e^u \geq 1$. Thus all the terms in $D_1D_{\bar 1}F$ can be bounded by the right hand side of (\ref{D1D1}), completing
the proof of (\ref{D1D1}).

%All terms of order $|\Delta u|$, $|D \bar{D} u|$, $|D D u|$ can be absorbed into the higher order term $|D D u|^2 + |D \bar{D} u|^2$. With this in mind, we can obtain the estimate
%\bea
%|D_1 D_{\bar{1}} F| &\leq& B (e^u + e^{-u}) \{ |D D u|^2 + |D \bar{D} u|^2 \} +B(e^u + e^{-u})(1 + |Du|^2) \{|D D u| + |D \bar{D} u|\} \nonumber\\
%&&+ B \{ e^u |D u| |D u_{\bar{1} 1}| + e^{-u} |D u_{\bar{1} 1}| \} + C\nonumber\\
%&\leq& B e^u  \left( |D D u|^2 + |D \bar{D} u|^2 + e^u|D u_{\bar{1} 1}| \right) +C.
%\eea
%In the last line, we used . Using the definition of $F$, we estimate
%\be \label{grad_F}
%|D F| \leq  B e^{u} |D u| \left( |D D u| + |D \bar{D} u| \right) + C. 
%\ee
Using the lower bound for $F$ in (\ref{F_lower}) and
the fact that $(e^u|Du|+e^{-u})^2\leq 2(e^{2u}|Du|^2+e^{-2u})\leq B e^{3u}$,
in view of the assumption $e^{-u} \leq 1$ and $|D u|^2 \leq e^u$ (\ref{C2_assumption}), we have
\begin{eqnarray*}
{|D F|^2 \over F} &\leq& {B e^{3u} \over F} \left( |D D u|^2 + |D \bar{D} u|^2 \right) + C\leq B e^u \left( |D D u|^2 + |D \bar{D} u|^2 \right) +C.
\end{eqnarray*}
Q.E.D.
\newline
\newline
{\it Proof of Lemma \ref{c2lemma2}.} This argument will adapt the proof of Proposition 9 in Guan-Ren-Wang \cite{GRW} to the complex setting. Recall that we are working at a point $p$ with $D G =0$, $g_{\bar{k} j}=\delta_{kj}$ and $u_{\bar{k} j}$, $g'$ diagonal, and $g'_{\bar{1} 1} \geq g'_{\bar{2} 2} \geq \dots \geq g'_{\bar{n} n}$.
%\be
%g'_{\bar{1} 1} \geq g'_{\bar{2} 2} \geq \dots \geq g'_{\bar{n} n}.
%\ee
The first step is the following computation, for $a \in \{1,2,\dots,n \}$ fixed.
\bea
& \ & -(1-\tau)F^{a \bar{a}} \left| {D_a g'_{\bar{1} 1} \over g'_{\bar{1} 1}} \right|^2  + {\phi' \over 2} F^{a \bar{a}} \left(|u_{a \bar{a}}|^2 +  \sum_k |u_{k a}|^2 \right)  \nonumber\\
&=& -(1-\tau)F^{a \bar{a}} \left| \phi' |D u|^2_a -2n\alpha  \psi' u_a \right|^2  + {\phi' \over 2} F^{a \bar{a}} \left(|u_{a \bar{a}}|^2 +  \sum_k |u_{k a}|^2 \right) \nonumber\\
&\geq& -2(1-\tau)F^{a \bar{a}} \left( \left| \phi' |D u|^2_a \right|^2 + \left|2n\alpha \psi' u_a \right|^2 \right) + {\phi' \over 2} F^{a \bar{a}} \left(|u_{a \bar{a}}|^2 +  \sum_k |u_{k a}|^2 \right) \nonumber\\
&\geq& F^{a \bar{a}} \bigg\{ \phi' \cdot \left( {1 \over 2} - 4(1-\tau)\phi'|D u|^2 \right) \left(|u_{a \bar{a}}|^2 +  \sum_k |u_{k a}|^2 \right) - C \bigg\}.
\eea
By our choice of $\phi$, we have (\ref{grad_upper}), (\ref{phi_upper}) and (\ref{phi_lower}), hence
\be
\phi' |D u|^2 \leq {1 \over 16}, \ \ \phi' > {1 \over C} > 0.
\ee
Thus we have for $a \in \{1,2,\dots,n \}$,
\be \label{lemma2_id}
-(1-\tau)F^{a \bar{a}} \left| {D_a g'_{\bar{1} 1} \over g'_{\bar{1} 1}} \right|^2  + {\phi' \over 2} F^{a \bar{a}} \left(|u_{a \bar{a}}|^2 +  \sum_k |u_{k a}|^2 \right) \geq  F^{a \bar{a}} \left( {1 \over 4} |u_{a \bar{a}}|^2 - C \right)
\ee
If $g'_{\bar{1} 1} \gg 1$, then
\be
u_{\bar{1} 1} = {1 \over 2n \alpha} \left(g'_{\bar{1} 1} - (e^u + f e^{-u})\right) \gg 1,
\ee
hence by letting $a=1$ in (\ref{lemma2_id})
\be
-(1-\tau)F^{1 \bar{1}} \left| {D_1 g'_{\bar{1} 1} \over g'_{\bar{1} 1}} \right|^2  + {\phi' \over 2} F^{1 \bar{1}} \left(|u_{1 \bar{1}}|^2 +  \sum_k |u_{k 1}|^2 \right) \geq 0
\ee
for $g'_{\bar{1} 1}$ sufficiently large. Thus to prove the lemma, it needs to be shown that
\be \label{lemma2_estimate}
{1-{\tau \over 2} \over g'_{\bar{1} 1}} \sum_{i >1 } |D_i g'_{\bar{1} 1}|^2 -(1-\tau) \sum_{i>1} F^{i \bar{i}} \left| {D_i g'_{\bar{1} 1} \over g'_{\bar{1} 1}} \right|^2  + \sum_{i>1} {\phi' \over 2} F^{i \bar{i}}  \left(|u_{i \bar{i}}|^2 +  \sum_k |u_{k i}|^2 \right)\geq 0.
\ee
To prove this estimate, we proceed by cases. Let 
\be
0< \varepsilon < {\tau \over 2 (1-\tau)}.
\ee
Case (A): $\sum_{i=2}^{n-1} g'_{\bar{i} i} < \varepsilon g'_{\bar{1} 1}$. In this case, we have $F^{n \bar{n}} = g'_{\bar{1} 1} + \sum_{i=2}^{n-1} g'_{\bar{i} i} < (1 + \varepsilon )g'_{\bar{1} 1}$.
%\be
%F^{n \bar{n}} = g'_{\bar{1} 1} + \sum_{i=2}^{n-1} g'_{\bar{i} i} < (1 + \varepsilon )g'_{\bar{1} 1}.
%\ee
Thus
\begin{eqnarray*}
&&{1-{\tau \over 2} \over g'_{\bar{1} 1}} \sum_{i >1 } |D_i g'_{\bar{1} 1}|^2 -(1-\tau) \sum_{i>1} F^{i \bar{i}} \left| {D_i g'_{\bar{1} 1} \over g'_{\bar{1} 1}} \right|^2\\ &\geq& {1- {\tau \over 2} \over g'_{\bar{1} 1}} \sum_{i >1 } |D_i g'_{\bar{1} 1}|^2 -(1-\tau) \sum_{i>1} F^{n \bar{n}} \left| {D_i g'_{\bar{1} 1} \over g'_{\bar{1} 1}} \right|^2 \nonumber \\
& \geq & \left(1 - {\tau \over 2}- (1-\tau)(1+\varepsilon) \right) {1 \over g'_{\bar{1} 1}} \sum_{i >1 } |D_i g'_{\bar{1} 1}|^2,\nonumber
\end{eqnarray*}
which is nonnegative by the choice of $\varepsilon$. This proves (\ref{lemma2_estimate}).
\newline
\newline
Case (B): $\sum_{i=2}^{n-1} g'_{\bar{i} i} \geq \varepsilon g'_{\bar{1} 1}$. In this case, we have $g'_{\bar{2} 2} \geq {\varepsilon \over n-2} g'_{\bar{1} 1}$.
%\be
%g'_{\bar{2} 2} \geq {\varepsilon \over n-2} g'_{\bar{1} 1}.
%\ee
For $g'_{\bar{1} 1}$ large enough,
\be
u_{\bar{2} 2} =  {1 \over 2n \alpha} \left(g'_{\bar{2} 2} - (e^u + f e^{-u}) \right) \geq {1 \over 2n \alpha}{\varepsilon \over n-2} g'_{\bar{1} 1} - C \geq {\varepsilon \over 4n(n-2) \alpha} g'_{\bar{1} 1}.
\ee
We divide case (B) into subcases.
\newline
\newline
Case (B1): $F^{2 \bar{2}} \geq 1$. By (\ref{lemma2_id})
\bea
&\ & -(1-\tau) \sum_{i>1} F^{i \bar{i}} \left| {D_i g'_{\bar{1} 1} \over g'_{\bar{1} 1}} \right|^2 + \sum_{i>1} {\phi' \over 2} F^{i \bar{i}}  \left(|u_{i \bar{i}}|^2 +  \sum_k |u_{k i}|^2 \right) \nonumber\\
& \geq & F^{2 \bar{2}} \left( {1 \over 4} |u_{2 \bar{2}}|^2 - C \right) + \sum_{i>2} F^{i \bar{i}} \left( {1 \over 4} |u_{i \bar{i}}|^2 - C \right) \nonumber\\
& \geq & \left( {1 \over 4} |u_{2 \bar{2}}|^2 - C \right) -C \sum_{i>2} F^{i \bar{i}}  \geq {\varepsilon^2 \over 4^3 (n (n-2) \alpha)^2}  (g'_{\bar{1} 1})^2  -C g'_{\bar{1} 1} -C \geq 0.
\eea
Case (B2): $F^{2 \bar{2}} < 1$. In this case,
\be
-g'_{\bar{n} n} > -1 + (g'_{\bar{1} 1} - g'_{\bar{2} 2}) + \sum_{i=2}^{n-1} g'_{\bar{i} i} \geq -1 + \varepsilon g'_{\bar{1} 1} \geq {\varepsilon \over 2} g'_{\bar{1} 1}.
\ee
\be
-u_{\bar{n} n} = {1 \over 2n \alpha} \left( -g'_{\bar{n} n} + (e^u + f e^{-u}) \right) \geq {\varepsilon \over 4n \alpha} g'_{\bar{1} 1}.
\ee
Note that the assumption of case (B) implies $F^{n \bar{n}} \geq (1+\varepsilon) g'_{1 \bar{1}}$.
%\be
%F^{n \bar{n}} \geq (1+\varepsilon) g'_{1 \bar{1}}.
%\ee
Another computation using (\ref{lemma2_id}) yields
\bea
& \ & -(1-\tau) \sum_{i>1} F^{i \bar{i}} \left| {D_i g'_{\bar{1} 1} \over g'_{\bar{1} 1}} \right|^2 + \sum_{i>1} {\phi' \over 2} F^{i \bar{i}}  \left(|u_{i \bar{i}}|^2 +  \sum_k |u_{k i}|^2 \right) \nonumber\\
&\geq& F^{n \bar{n}} \left( {1 \over 4} |u_{n \bar{n}}|^2 - C \right) + \sum_{i=2}^{n-1} F^{i \bar{i}} \left( {1 \over 4} |u_{i \bar{i}}|^2 - C \right) \nonumber\\
&\geq& F^{n \bar{n}} \left( {\varepsilon^2 \over 4^3(n \alpha)^2} |g'_{1 \bar{1}}|^2 - C \right) - C g'_{1 \bar{1}}\nonumber\\
&\geq& (1+\epsilon) g'_{1 \bar{1}} \left( {\varepsilon^2 \over 4^3(n \alpha)^2} |g'_{1 \bar{1}}|^2 - C \right) - C g'_{1 \bar{1}} \geq 0.
\eea
This establishes (\ref{lemma2_estimate}), and thus proves Lemma \ref{c2lemma2}. Q.E.D.

\section{The $C^{2,\eta}$ Estimate}
\setcounter{equation}{0}

At this point, we have shown the {\it a priori} $C^2$ estimates (\ref{C2estimate1}) for equation (\ref{scalar_eqn}), under the assumption of a sharp $C^1$ upper bound (\ref{C2_assumption}). This $C^2$ estimate implies that the equation is uniformly elliptic and that it is also a concave operator. We would like to apply the Evans-Krylov theorem \cite{Evans, Krylov, Siu} to show the $C^{2, \eta}$ bound. However, we cannot apply the standard theorem directly.

\par
In fact, equation (\ref{scalar_eqn}) is of the following form
\begin{eqnarray*}
\sigma_2\left(\chi_{\bar{j} k}(z, u) + u_{\bar{j} k}\right) = \varphi(z, u, Du).
\end{eqnarray*}
By the {\it a priori} $C^2$ estimate, we have uniform bounds for the complex Hessian $\partial\bar\partial u$ and hence for $\Delta u$. This implies that $u\in C^{1, \theta}$ for some $\theta\in (0, 1)$. Therefore, the function $\varphi(z, u, Du)=$ right hand side of equation (\ref{scalar_eqn}) is only $C^{\theta}$ even if $f$ and $\mu$ are smooth on $X$. Thus, the standard Evans-Krylov theorem is not directly applicable as it requires a $C^{1, 1}$ bound for $\varphi$, which depends on the $C^3$ norm of $u$ in our case.

\par
The $C^{2, \eta}$ regularity for the complex Monge-Amp\`ere equations with only H\"older continuous right hand side was obtained by Dinew-Zhang-Zhang \cite{DZZ} for $u\in C^{1, 1}$. The assumption on $u$ was weaken to be $\Delta u\in L^{\infty}$ by Wang \cite{Wang} and it was later extended to more general settings by Tosatti-Wang-Weinkove-Yang \cite{TWWY}. Indeed, our setup here fits well into the general picture in \cite{TWWY} (Theorem 1.1 for equation (1.4) in \cite{TWWY}). We note that our $\chi_{\bar{j} k} = \left(e^{u}+fe^{-u}\right) g_{\bar{j} k} \in C^{1, \theta}$ and $\varphi \in C^{\theta}$. And thus we can apply their main result to conclude the following $C^{2, \eta}$ bound for $u$. We refer the reader to \cite{TWWY} for details.

\begin{theorem} \label{C2alphaestimate}
Let $u\in C^2(X)$ be a solution to (\ref{scalar_eqn}) with normalization condition  (\ref{A_normalize}). Then, there exist positive constants $0<\eta<1$ and $C$ depending on $n$, $(X, g)$, $\| u \|_{L^\infty}$, $\|D u\|_{L^\infty}$, $\|\Delta u\|_{L^\infty}$, $\|f\|_{C^3}$, $\|\mu\|_{C^1}$ and $\alpha$ such that
\be
\|u\|_{C^{2, \eta}(X)} \leq C.
\ee
\end{theorem}

\bigskip

\section{Non-Degeneracy and Sharp Gradient Bounds} \label{lower_b_section}
\setcounter{equation}{0}
\par
In order to solve equation (\ref{scalar_eqn}) subject to normalization condition (\ref{A_normalize}), one can use the method of continuity. This can be done by introducing the parameter $t$, and replacing $f$ by $tf$ and $\mu$ by $t \mu$. 
\bea \label{method_cont}
e^{-2u} F
&=&
\kappa_c \{1-4\alpha e^{-u} |Du|^2 \} + 4 \alpha \kappa_c \bigg \{  tf e^{-3u} |D u|^2  - t e^{-3u} (g^{i \bar{k}}f_i u_{\bar{k}} + g^{i \bar{k}} f_{\bar{k}} u_i) \bigg \}
\nonumber\\
&&
+ \kappa_c te^{-2u}  \bigg\{ 2f + f^2 e^{-2u} + 4 \alpha e^{-u} \Delta f \bigg\}
-2n\alpha t e^{-2u} \mu.
\eea
We see that when $t=0$, the equation admits the trivial solution $u=- \log A$, and the right hand side is equal to $\kappa_c$. The issue addressed in this section is whether the right-hand side can degenerate to zero as $t$ tends to $t=1$. For simplicity, we shall suppress the parameter $t$ in our computations and write $f$ instead of $tf$ and $\mu$ instead of $t \mu$. The theorem of Fu-Yau \cite{FY} is the following.
\begin{theorem}[Fu-Yau \cite{FY}] \label{lower_bdd_thm}
Let the dimension of $X$ be equal to $n=2$. For any $\delta >0$, there exists $A_\delta>0$ depending on $(X,\o)$, $f$ $\alpha$, $\mu$, such that if $A < A_\delta$, then for any solution $u$ of the Fu-Yau equation (\ref{scalar_eqn}) with normalized condition (\ref{A_normalize}), there holds
\be
e^{-2u} F \geq \kappa_c -\delta.
\ee

\end{theorem}

\medskip
In the rest of this section, we investigate the non-degeneracy estimate for the higher dimensional case. As mentioned in the Introduction,
we follow the idea of Fu-Yau closely, but we work with general coordinate systems rather than the adapted ones with $\nabla u=(u_1,0,\cdots,0)$ used by Fu and Yau. This allows us a simplified and more transparent derivation of the Fu-Yau results for $n=2$, and a clearer picture of why their arguments are not strong enough for higher dimensions. Following Fu-Yau, we apply the maximum principle to the following function
\bea
G=1-4\alpha e^{-u}|D u|^2+4\alpha e^{-\e u}-4\alpha e^{-\e\,{\rm inf} \,u}.
\eea

\subsection{First computation of $F^{j\bar k}D_jD_{\bar k}G$}

We begin by computing $F^{j\bar k}D_jD_{\bar k}(-4\alpha e^{-u}|D u|^2)$.
Because we shall ultimately evaluate this expression as a critical point of $G$, where 
\bea
\label{critical}
D(e^{-u}|D u|^2)=D(e^{-\e u})
\eea
it is advantageous to express $D_jD_{\bar k}(-4\alpha e^{-u}|D u|^2)$ in terms of $D(e^{-u}|D u|^2)$ as much as possible. Thus we write
\bea
F^{j\bar k}D_jD_{\bar k}
(-4\alpha e^{-u}|D u|^2)
&=&
4\alpha F^{j\bar k}D_{\bar k}u D_j(e^{-u}|D u|^2)
+
4\alpha F^{j\bar k}D_j uD_{\bar k}(e^{-u}|D u|^2)
\nonumber\\
&&
+
4\alpha  e^{-u}|D u|^2\,F^{j\bar k} D_ju D_{\bar k}u
-
4\alpha e^{-u}F^{j\bar k}D_jD_{\bar k}|D u|^2
\nonumber\\
&&
+4\alpha (F^{j\bar k}D_jD_{\bar k}u)\,e^{-u}|D u|^2.
\eea
On the other hand, a straightforward computation gives
\bea
F^{j\bar k}D_j D_{\bar k}(4\alpha e^{-\e u})
=
4\alpha \e^2 |D u|_F^2 e^{-\e u}
-
4\alpha\e (F^{j\bar k}D_jD_{\bar k}u) e^{-\e u}.
\eea
and thus
\bea
F^{j\bar k}D_jD_{\bar k}G
&=&
4\alpha F^{j\bar k}D_{\bar k}u D_j(e^{-u}|D u|^2)
+
4\alpha F^{j\bar k}D_j uD_{\bar k}(e^{-u}|D u|^2)
\nonumber\\
&&
+
4\alpha  (e^{-u}|D u|^2+\e^2 e^{-\e u})\,|Du|_F^2
-
4\alpha e^{-u}F^{j\bar k}D_jD_{\bar k}|D u|^2
\nonumber\\
&&
+4\alpha (F^{j\bar k}D_jD_{\bar k}u)\,(e^{-u}|D u|^2-\e e^{-\e u})
\eea
where we have introduced the notation $|Du|_F^2=F^{j\bar k}D_j u D_{\bar k}u$.
%\bea
%|Du|_F^2=F^{j\bar k}D_j u D_{\bar k}u.
%\eea
\par
We can now substitute in the critical point equation (\ref{critical}) of $G$, and obtain
\bea
F^{j\bar k}D_jD_{\bar k}G
&=&
4\alpha (e^{-u}|Du|^2-2\e e^{-\e u}+\e^2 e^{-\e u})|Du|_F^2
\nonumber\\
&&
+ 4 \alpha (F^{j\bar k}D_jD_{\bar k}u) \,( e^{-u}|D u|^2- \e e^{-\e u})
-
4\alpha e^{-u}F^{j\bar k}D_jD_{\bar k}|D u|^2
\eea
Both expressions $F^{j\bar k}D_jD_{\bar k}u$ and $F^{j\bar k}D_jD_{\bar k}|Du|^2$
have been computed in section \S 2 and are found in equations (\ref{FDDu}) and (\ref{DD|Du|^2}). Substituting in the formulas derived there, we obtain
\bea
\label{G1}
F^{j\bar k}D_jD_{\bar k}G
&=&
4\alpha (e^{-u}|Du|^2-2\e e^{-\e u}+\e^2 e^{-\e u})|Du|_F^2
-4\alpha e^{-u}(|DDu|_{Fg}^2+|D\bar Du|_{Fg}^2)
\nonumber\\
&&
+
\bigg\{ {4\over n}F
-
{2(n-1)\over n}(e^u+fe^{-u}){\rm Tr}\,h \bigg\} (e^{-u}|D u|^2-\e e^{-\e u})
\nonumber\\
&&
+
{2(n-1)\over n}{\rm Tr}\,h\,e^{-u}g^{\ell\bar m}
\big\{\p_\ell(e^u+fe^{-u})\p_{\bar m}u+\p_{\bar m}(e^u+fe^{-u})\p_\ell u\big\}
\nonumber\\
&&
-{2\over n} e^{-u}g^{\ell\bar m}\big\{\p_\ell F \p_{\bar m}u+\p_{\bar m} F \p_\ell u\big\}
+4\alpha e^{-u}\tilde g_{\bar\ell m}R^{m\bar\ell p\bar q}\p_pu \p_{\bar q}u.
\eea

%\subsection{Blocki's partial cancellation between $|DD u|_{Fg}^2$ and
%$|Du|^2|Du|_{Fg}^2$}

We now make use of a key partial cancellation, observed by Blocki in his proof of $C^1$ estimates for the Monge-Amp\`ere equation \cite{Bl} (see also \cite{Guan, Zhang}, and \cite{PS, PSS} for other applications of this partial cancellation), between $|DDu|_{Fg}^2$ and $|Du|^2|Du|_{Fg}^2$, which is the following. At a critical point of $G$, the relation (\ref{critical}) implies
\bea
\label{critical1}
g^{i\bar j}D_pD_i u D_{\bar j}u=
-g^{i\bar j}D_{\bar j}D_pu D_iu
+(|Du|^2-\e e^{(1-\e)u})D_pu.
\eea
We can now estimate $|DDu|_{Fg}^2$ from below by
\bea
|DDu|_{Fg}^2&\geq& {1\over |Du|^2}|g^{i\bar j}D_p D_i u D_{\bar j}u|_F^2
={1\over |Du|^2} \bigg| g^{i\bar j}D_{\bar j}D_pu D_iu
-(|Du|^2-\e e^{(1-\e)u})D_pu \bigg|_F^2
\nonumber\\
&=&
|Du|^2|Du|_F^2 +\e^2 e^{2(1-\e)u}{|Du|_F^2\over |Du|^2}-2\e |Du|_F^2 e^{(1-\e)u}
\nonumber\\
&&
+
{1\over |Du|^2}|g^{i\bar j}D_iu D_{\bar j}D_p u|_F^2
-
{2\over |Du|^2}
(|Du|^2-\e e^{(1-\e)u}){\rm Re}(F^{p\bar q}g^{i\bar j}D_{\bar q}uD_iu D_{\bar j}D_pu).
\nonumber
\eea
The terms $|Du|^2|Du|_F^2$ and
$-2\e e^{(1-\e)u}|Du|_F^2$ will cancel out similar terms in $F^{j\bar k}D_jD_{\bar k}G$. The expression ${\rm Re}(F^{p\bar q}g^{i\bar j}D_{\bar q} u D_iu D_{\bar j}D_pu)$ can be rewritten as
\bea
{\rm Re}(F^{p\bar q}g^{i\bar j}D_{\bar q}uD_iu D_{\bar j}D_pu)
&=&
{1\over 2n\alpha}
{\rm Re}(F^{p\bar q}g^{i\bar j}D_{\bar q}uD_iu g_{\bar jp}')
-
{e^u+fe^{-u} \over 2n\alpha}|Du|_F^2
\nonumber\\
&=&
{1\over 2n\alpha}\{F|Du|^2-\sum_{j=1}^n\sigma_2(\lambda'|j)|u_j|^2-(e^u+fe^{-u})|Du|_F^2\} \nonumber
\eea
by going to coordinates where $g_{\bar kj}=\delta_{\bar kj}$, and $u_{\bar kj}$ is diagonal at the point $p$ where the function $G$ attains its minimum. Here $\sigma_k(\lambda'|j)$ denotes the $k$-th symmetric function of the $(n-1)\times(n-1)$ diagonal matrix with eigenvalues $\lambda'_m$, $m\not=j$.

\smallskip

We also make use of the other term $|D\bar D u|_{Fg}^2$, which we rewrite as
\begin{eqnarray*}
|D\bar D u|_{Fg}^2&=&{F^{p\bar q}g^{i\bar j}\over (2n\alpha)^2}
(g_{\bar jp}'-(e^u+fe^{-u})g_{\bar jp})(g_{\bar qi}'-(e^u+fe^{-u})g_{\bar qi})
\nonumber\\&=&{1\over (2n\alpha)^2}
(|g'|_{Fg}^2 -4(e^u+fe^{-u}) F+(n-1)(e^u+fe^{-u})^2{\rm Tr}\,h).
\end{eqnarray*}
Again using the above coordinates, we can work out a more explicit expression for $|g'|_{Fg}^2$,
\bea
-{e^{-u}\over n^2\alpha}
|g'|_{Fg}^2
&=&
-{e^{-u}\over n^2\alpha}\sum_{j=1}^n\tilde\lambda_j(\lambda_j')^2
=
-{e^{-u}\over n^2\alpha}(F\sum_{j=1}^n\lambda_j'-\sum_{j=1}^n\sigma_2(\lambda'|j)\lambda_j')
\nonumber\\
&=&
-{e^{-u}\over n^2\alpha}
(F{\rm Tr}\,h-\sum_{j=1}^n (\sigma_3(g')-\sigma_3(\lambda'|j))
=
-{e^{-u}\over n^2\alpha}F\,{\rm Tr}\,h
+{e^{-u}\over n^2\alpha}3\sigma_3(g').
\nonumber
\eea
Thus we find, at a critical point of the test function $G$,
\bea \label{before_using_eqn}
F^{j\bar k}D_jD_{\bar k}G
&\leq&
\bigg\{-{4\over n}e^{-2u}(e^u+fe^{-u})|Du|^2+{4\over n}\e e^{-(1+\e) u}(e^u+fe^{-u})
-4\alpha\e^2 e^{-2\e u} \nonumber\\
&&+4\alpha\e^2
e^{-\e u}e^{-u}|Du|^2\bigg\} \, e^u{|Du|_F^2\over |Du|^2}-
{4\alpha e^{-u}\over |Du|^2}|g^{i\bar j}D_iu D_{\bar j}D_p u|_F^2
\nonumber\\
&&
+(e^{-u}|Du|^2-\e e^{-\e u})\bigg\{
{8\over n}F-{4\over n}\sum_{j=1}^n\sigma_2(\lambda'|j){|u_j|^2\over|Du|^2}
-
2{n-1\over n}(e^u+fe^{-u}){\rm Tr}\,h\bigg\}
\nonumber\\
&&
-{e^{-u}\over n^2\alpha}F\,{\rm Tr}\,h
+{e^{-u}\over n^2\alpha}3\sigma_3(g')
+{4\over n^2\alpha}e^{-u}(e^u+f e^{-u})F
-
{n-1\over n^2\alpha}e^{-u}(e^u+fe^{-u})^2{\rm Tr}\,h
\nonumber\\
&&
+
{2(n-1)\over n}{\rm Tr}\,h\,e^{-u} \, 2 \Re \, \<D(e^u + f e^{-u}),Du\> 
-{2\over n}e^{-u} \,
2{\rm Re}\<DF,Du\> \nonumber\\
&&
+4\alpha e^{-u}\tilde g_{\bar\ell m}R^{m\bar\ell p\bar q}\p_pu \p_{\bar q}u.
\eea

\subsection{Using the equation}

So far, we have not used the equation (\ref{scalar_eqn}). We shall now use it to
evaluate and simplify the preceding estimate for
$F^{j\bar k}D_jD_{\bar k}G$. It is convenient to rewrite the equation (\ref{scalar_eqn}) in the following form
\bea
\label{scalar_eqn1}
F
=
{n(n-1)\over 2}e^{2u}(1-4\alpha e^{-u}|Du|^2)-2n\alpha \nu
\eea
where the function $\nu$ is defined to be
\bea
\nu
&=&
\mu
-(n-1)fe^{-u}|Du|^2-{n-1\over 2\alpha}f-{n-1\over 4\alpha}e^{-2u}f^2
\nonumber\\
&&
-
(n-1)e^{-u}(\Delta f-g^{j\bar k}(D_jf D_{\bar k}u+D_{\bar k}uD_jf)).
\eea
At a critical point (\ref{critical}) for $G$, we have
\bea
\p_{\bar \ell} F
&=&
n(n-1)e^{2u}\p_{\bar\ell}u(1-4\alpha e^{-u}|Du|^2)
+
{n(n-1)\over 2}e^{2u}\p_{\bar\ell} (-4\alpha e^{-u}|Du|^2)-2n\alpha\p_{\bar\ell}\nu
\nonumber\\
&=&
2F\,\p_{\bar\ell}u-4n\alpha \p_{\bar\ell}u\nu
+
2\alpha{n(n-1)}\e \p_{\bar\ell}u e^{(2-\e) u}
-2n\alpha\p_{\bar\ell}\nu.
\eea
The preceding expression is unwieldy if we write it down in full. To avoid unnecessary details, it is convenient to introduce the following groups of expressions:

\medskip
$\bullet$ The group ${\cal E}_0$ consists of the following expressions
\bea
e^{-u}{\rm Tr}\,h,
\quad
e^{-2u} F,
\quad
e^{-3u}\sigma_3(g'),
\quad
e^{-u}{|Du|_F^2\over|Du|^2},
\quad
e^{-2u} \sum_j\sigma_2(\lambda'|j){|u_j|^2\over |Du|^2},
\quad 1
\eea
where $\sigma_3$ and $\sigma_2(\lambda'|j)$ denote the symmetric functions of the eigenvalues of the matrix $g_{\bar kj}'$.

\smallskip
$\bullet$ The group ${\cal E}_1$ consists of expressions of the form
\bea
\e e^{-\e u}\,\Phi,
\qquad
\Phi\in {\cal E}_0.
\eea

\smallskip
$\bullet$ The group ${\cal E}_2$ consists of expressions of the form
\bea
c\,\Phi,
\qquad
\textit{ with } \Phi\in {\cal E}_0
\ \textit{ and }\
c<<\e e^{-\e u}
\eea
where the inequality indicated on the coefficient $c$ should hold for $\e<<1$ and $A<<1$. For example, any function of the form $e^{-2u}v$ with $v$ a bounded function can be classified into the group ${\cal E}_2$. Another example is the expression $4\alpha\e^2 e^{-\e u}|Du|_F^2$, which can be viewed as belonging to $e^{2u} {\cal E}_2$, since
\bea
4\alpha\e^2 e^{-\e u}|Du|_F^2
=
e^{2u} \bigg\{ \e^2 e^{-\e u}(4\alpha e^{-u}|Du|^2)\, e^{-u}{|Du|_F^2\over |Du|^2} \bigg\}
\eea
and the expression $4\alpha e^{-u}|Du|^2$ is bounded as we vary $A$ by the $C^1$ estimate (\ref{gradient_bound}).

We now claim that
\bea \label{using_eqn}
-{2\over n}e^{-u}
2{\rm Re}\<DF,Du\>
=-{8\over n}F\, e^{-u} |Du|^2
-8\alpha(n-1)\e e^{(2-\e) u}(e^{-u}|Du|^2)
\eea
modulo terms of the form $e^{2u} {\cal E}_2$. Indeed, absorbing all terms $e^{-u}|Du|^2$ into $O(1)$ yields
\bea
{-2 \over n} e^{-u} (-2n \alpha) 2 \Re \< D\nu, D u \> &=& - 4 \alpha (n-1) fe^{-u} \, 2 \Re \<Du,D(e^{-u}|Du|^2) \> \nonumber\\
&& + 4 \alpha (n-1) e^{-u} \, 2 \Re \< Df, D(e^{-u} |Du|^2) \> + O(1). \nonumber
\eea
Applying the critical point equation (\ref{critical}), we see that this term is of the form $e^{2u} {\cal E}_2$.

%This is easily verified by differentiating $\p\nu$ explicitly. Note that this produces terms of the form $g^{m\bar k}D_jD_{\bar k}u D_mu$, which can be rewritten in terms of $g_{\bar kj}'$, and estimated by ${\rm Tr}\,h$, as well as terms of the form $g^{\bar km}D_jD_mu D_{\bar k}u$. For these, we use the critical point equation (\ref{critical1}), to rewrite them in terms of $D_jD_{\bar k}u$, and ultimately in terms of $g_{\bar kj}'$ as well. 

\subsubsection{The expression for $F^{j\bar k}D_jD_{\bar k}G$ up to ${\cal E}_2$ terms}

It is now easy to clean up considerably the expression for $F^{j\bar k}D_jD_{\bar k}G$. Up to $e^{2u} {\cal E}_1$ and $e^{2u} {\cal E}_2$ terms, (\ref{before_using_eqn}) is
\bea
F^{j\bar k}D_jD_{\bar k}G
&\leq&
- \{ {4\over n}e^{-u}|Du|^2 \} e^u {|Du|_F^2\over |Du|^2}
+
e^{-u}|Du|^2\bigg\{{8\over n}F-{4\over n}\sum_j\sigma_2(\lambda'|j){|u_j|^2\over |Du|^2}
-
2{n-1\over n}e^u{\rm Tr}\,h\bigg\}
\nonumber\\
&&
-
{4\alpha e^{-u}\over|Du|^2}|g^{i\bar j}D_iuD_{\bar j}D_pu|_F^2
-
{e^{-u}\over n^2\alpha}F\,{\rm Tr}\,h
+
{e^{-u}\over n^2\alpha}3\sigma_3
+
{4\over n^2\alpha}F-{n-1\over n^2\alpha}e^u{\rm Tr}\,h
\nonumber\\
&&
+
4{n-1\over n} {\rm Tr}\,h |Du|^2
-
{2\over n} e^{-u} \, 2{\rm Re}\<DF,Du\>.
\nonumber
\eea
We can now make use of the formula (\ref{using_eqn}) for $\<D F,Du\>$ modulo $e^{2u} {\cal E}_1$ and $e^{2u} {\cal E}_2$ obtained in the previous section. The expression $8F/n$ in the top line cancels out. Regrouping terms in terms of $|Du|_F^2/|Du|^2$, ${\rm Tr}\,h$ and $F$, we have
\bea
F^{j\bar k}D_jD_{\bar k}G
&\leq&
e^u{\rm Tr}\,h
\bigg\{2{n-1\over n}e^{-u}|Du|^2 -{n-1\over n^2\alpha}-{e^{-2u}F\over n^2\alpha}\bigg\}
+
 {4\over n^2\alpha}F -\{{4\over n}e^{-u}|Du|^2 \} e^u{|Du|_F^2\over |Du|^2}
\nonumber\\
&&
-{4 e^{-u} |D u|^2 \over n}\sum_j\sigma_2(\lambda'|j){|u_j|^2\over |Du|^2}
+{e^{-u}\over n^2\alpha}3\sigma_3
-{4\alpha e^{-u}\over |Du|^2}
|g^{i\bar j}D_iu D_{\bar j}D_pu|_F^2.
\eea
We can now eliminate systematically $4\alpha e^{-u}|Du|^2$ using the equation
\bea
4\alpha e^{-u}|Du|^2
=1-{e^{-2u}F\over\kappa_c}\quad {\rm modulo}\ {\cal E}_2.
\eea
where it is convenient to introduce the critical value $\kappa_c$ as in (\ref{definekc}).
The coefficient of $e^u {\rm Tr}\,h$ above becomes
\bea
2{n-1\over n} e^{-u}|Du|^2 -{n-1\over n^2\alpha}-{e^{-2u}F\over n^2\alpha}
=
{n-1\over n\alpha}
\big\{{1\over 2}-{1\over n}-{e^{-2u}F\over\kappa_c}\big\}.
\eea
We multiply by $e^{-2u}$ and summarize the previous calculations in the following inequality, 
\bea
\label{inequality}
(F^{j\bar k}D_jD_{\bar k}G)e^{-2u}
&\leq&
\bigg\{{1\over 2}-{1\over n}-{e^{-2u}F\over\kappa_c}\bigg\}{n-1\over n\alpha}e^{-u}\,{\rm Tr}\,h
+
{3\over n^2\alpha}e^{-3u}\sigma_3+{4\over n^2\alpha}e^{-2u}F
\nonumber\\
&&
-\bigg\{ 1-{e^{-2u}F \over\kappa_c} \bigg\} {1\over n\alpha}e^{-u}{|Du|_F^2\over|Du|^2}
-{1\over n\alpha} \bigg\{ 1-{e^{-2u}F \over\kappa_c} \bigg\}e^{-2u}\sum_{j=1}^n\sigma_2(\lambda'|j){|u_j|^2\over|Du|^2}
\nonumber\\
&&
-{4\alpha e^{-3u}\over|Du|^2}|g^{i\bar j}D_iu D_{\bar j} D_p u|_F^2
\eea
modulo terms in groups ${\cal E}_1$ and ${\cal E}_2$. The terms in group ${\cal E}_1$ in the expression for $(F^{j\bar k}D_jD_{\bar k}G)e^{-2u}$ come from (\ref{before_using_eqn}) and (\ref{using_eqn}) and can be worked out to be
\bea
&&
\e e^{-2u} e^{-\e u}
\bigg\{
2{n-1\over n}e^u\,{\rm Tr}\,h
-
{8\over n}F+
{4\over n}e^u {|Du|_F^2\over |Du|^2}+{4\over n}\sum_j \sigma_2(\lambda'|j){|u_j|^2\over |Du|^2}
-2(n-1)e^{2u}(1-{e^{-2u}F\over \kappa_c})\bigg\}
\nonumber\\
&&
=
\e e^{-2u} e^{-\e u}
\bigg\{
2{n-1\over n}e^u\,{\rm Tr}\,h
-
{4\over n}F-2(n-1)e^{2u}
+
{4\over n}e^u {|Du|_F^2\over |Du|^2}
+{4\over n}\sum_j \sigma_2(\lambda'|j){|u_j|^2\over |Du|^2}\bigg\}.
\eea
An explicit expression for the expression $|g^{i\bar j}u_iu_{\bar jp}|_F^2$
occurring above is
\bea
-{4 \alpha e^{-3u} \over |D u|^2}\left|  g^{i \bar{j}} u_{i} u_{\bar{j} p} \right|^2_F &=& - {e^{-3u} \over n^2 \alpha |D u|^2} \left| g^{i \bar{j}} u_i g'_{\bar{j} p} \right|_F^2 - {e^{-3u} (e^u+fe^{-u})^2 |D u|^2_F \over |D u|^2 n^2 \alpha} \nonumber\\
&&+ {2e^{-3u}(e^u+fe^{-u}) \over |D u|^2 n^2 \alpha}F^{p \bar{q}} g^{i \bar{j}} u_{\bar{q}}u_i g'_{\bar{j} p} \nonumber\\
&=&
-{e^{-3u}\over n^2\alpha}F\,{\rm Tr}\,h
+
{2 e^{-2u} \over n^2\alpha}F+
({e^{-2u}F\over n^2\alpha}-{1\over n^2\alpha})e^{-u}{|Du|_F^2\over|Du|^2}
\nonumber\\
&&
+
{e^{-3u}\over n^2\alpha}\sigma_3
-
{e^{-3u}\over n^2\alpha}
\sum_{j=1}^n\sigma_3(\lambda'|j){|u_j|^2\over|Du|^2}
-
{2 e^{-2u} \over n^2\alpha}\sum_{j=1}^n
\sigma_2(\lambda'|j){|u_j|^2\over|Du|^2}
\nonumber
\eea
modulo terms in group ${\cal E}_2$. Here we used
\bea
\left| g^{i \bar{j}} u_i g'_{\bar{j} p} \right|_F^2 = \sum_{j=1}^n \tilde{\lambda}_j \lambda'_j |u_j|^2 \lambda'_j =   \sum_{j=1}^n (F-\sigma_2(\lambda'|j)) |u_j|^2 ({\rm Tr}\,h -\tilde{\lambda_j}). \nonumber
\eea
Combining this expression with the previous two expressions, we obtain the following

\begin{theorem}
\label{minimum}
Let $p\in X$ be a point where the function $G$ achieves its minimum. Set
\bea 
\kappa_p=(e^{-2u}F)(p),
\qquad
\theta=2\alpha \e e^{-\e u}(p).
\eea
Then we have
\bea \label{minimum_ineq}
0
&\leq&
\bigg\{{1\over 2}-{1\over n}-{3\over 2}{\kappa_p\over\kappa_c}+\theta\bigg\}{n-1\over n}e^{-u}\,{\rm Tr}\,h
\nonumber\\
&&+\bigg\{{n+1\over n}(\kappa_p{1\over n-1}-1)+2\theta\bigg\}{1\over n}e^{-u}{|Du|_F^2\over |Du|^2}
\nonumber\\
&&
+
\bigg\{ {6\over n^2}-{2\over n}\theta \bigg\} \kappa_p
-
(n-1)\theta
-
{e^{-3u}\over n^2}
\sum_{j=1}^n\sigma_3(\lambda'|j){|u_j|^2\over|Du|^2}
\nonumber\\
&&
+
{4\over n^2}e^{-3u}\sigma_3
-{e^{-2u} \over n} \bigg\{ {n+2\over n}-{\kappa_p\over\kappa_c}-2\theta \bigg\}
\sum_{j=1}^n\sigma_2(\lambda'|j){|u_j|^2\over|Du|^2}
\eea
up to terms in group ${\cal E}_2$.
\end{theorem}

\subsection{A simplified Fu-Yau argument in dimension $n=2$}

We can now rederive the following key estimate of Fu-Yau \cite{FY} when $n=2$ (and hence $\kappa_c=1$): for any $\delta>0$, there exists
$A_\delta>0$ so that, if $A<A_\delta$, then the minimum $\kappa={\rm min}_X(e^{-2u}F)$ at any time $t$
satisfies the lower bound
\bea
\label{kappa}
\kappa>1-2\delta.
\eea
\par
Indeed, fix $\delta>0$, with $\delta<<1$. Recall that the test function $G(z)$ assumes its minimum at a point $p$, and set
$\kappa_p=(e^{-2u}F)(p)$. In view of the $C^0$ estimate, 
\begin{eqnarray*}
e^{-2u}F= \kappa_c G+O(A^{\e}),
\end{eqnarray*}
 and hence
$\kappa \geq \kappa_p+O(A^\e)$. Thus it suffices to show that (\ref{kappa}) holds with $\kappa$ replaced by $\kappa_p$.
It also suffices to show that if $\kappa_p>1/4$, then $\kappa_p> 1-\delta$ for $A_\delta$ small enough.
This is because $\kappa_p=\kappa=1$ when $t=0$, as discussed in (\ref{method_cont}). As $t$ varies,
$\kappa$ cannot reach $1/2$, since the first time it does so, we would have then $\kappa_p>1/4$ (for $A_0$ small enough),
and hence $\kappa>1-3\delta$, which is a contradiction. But then $\kappa>1/2$ for all time, and hence $\kappa>1-2\delta$
for all time, as desired.

\smallskip
We now argue by contradiction. Assume that $\kappa_p>1/4$. If $\kappa_p>1-\delta/2$, we are done, so we assume that
$\kappa_p\leq 1-\delta/2$.  In dimension $n=2$, $\sigma_3$ and $\sigma_2(\lambda'|j)$ all vanish. Incorporating the error terms in ${\cal E}_1$ and ${\cal E}_2$,
the inequality (\ref{minimum_ineq}) implies, for $A$ small enough,
\bea
c_1 e^{-u}{\rm Tr}\,h(p)+c_2e^{-u}{|Du|_F^2\over |Du|^2} (p)\leq c_3 \kappa_p +c_4
\eea
where $c_1$, $c_2$, $c_3$, $c_4$ are strictly positive constants, depending only on $\delta$. Since $\kappa_p$ is bounded by an absolute constant,
it follows that 
\begin{eqnarray*}
\label{gradient1}
e^{-u}\left({\rm Tr}\,h(p)+{|Du|_F^2\over |Du|^2}(p)\right)\leq c_5,
\end{eqnarray*}
 where $c_5$ is a constant depending only on $\delta$. This implies that all terms in $\theta^{-1} {\cal E}_2$ can be bounded by $c$, where $c$ is a constant that can be made arbitrarily small by taking $\e$ and $A$ to be small.

\smallskip
Going back again to the inequality (\ref{minimum_ineq}), we can bound the term $|Du|_F^2$ as follows,
\bea
\bigg\{{n+1\over n}(\kappa_p{1\over n-1}-1)+2\theta\bigg\}{1\over n\alpha}e^u{|Du|_F^2\over |Du|^2}
\leq
\bigg\{{n+1\over n}(\kappa_p{1\over n-1}-1)+2\theta\bigg\}{1\over n\alpha}e^u
\lambda_1'
\eea
where $\lambda_1'$ is either the largest or the lowest eigenvalue of $g_{\bar pq}'$, depending on the sign of the coefficient. In dimension $n=2$, $\kappa_c=1$, and Theorem \ref{minimum} implies, modulo additive terms of order ${\cal E}_2$,
\bea
0
&\leq& (-{3\over 4}\kappa_p+{1\over 2}\theta)(\lambda_1'+\lambda_2')e^{-u}
+
({3\over 4}(\kappa_p-1)+\theta)e^{-u}\lambda_1'
+
({3\over 2}-\theta)\kappa_p-\theta 
\nonumber\\
&=&
-({3\over 4}-3{\theta\over 2})e^{-u}\lambda_1'
-({3\over 4}\kappa_p-{\theta\over 2})e^{-u}\lambda_2'
+
({3\over 2}-\theta)\kappa_p-\theta.
\eea
Since $a_1\lambda_1'+a_2\lambda_2'\geq 2\sqrt {a_1a_2}\,\sqrt{\lambda_1'\lambda_2'}$ for any $a_1,a_2\geq 0$, and since $\lambda_1'\lambda_2'=e^{2u}\kappa_p$, we obtain
\bea
\label{inequalitykp}
({3\over 4}-3{\theta\over 2})^{1\over 2}
({3\over 4}\kappa_p-{\theta\over 2})^{1\over 2}
\kappa_p^{1\over 2}
\leq
{3\over 4}\kappa_p-{\theta\over 2}(\kappa_p+1).
\eea
The leading term $\kappa_p^2$ cancels upon squaring both sides. Since we have assumed that $\kappa$ is bounded  away from $0$, we can also divide by $\kappa_p$ and the error terms of type ${\cal E}_2$ will remain of type ${\cal E}_2$. We obtain, discarding terms of order $\theta^2$ and dividing through by $\theta\kappa_p$, $-\kappa_p-{1\over 3}\leq -{2\over 3}(\kappa_p+1)$, or equivalently,
\bea
\kappa_p\geq 1
\eea
modulo additive constants which can be made arbitrarily small by taking $A$ small. This establishes the desired lower bound for $\kappa$. 

\

To finish the discussion on dimension $n=2$, we note that Theorem \ref{lower_bdd_thm} implies the sharp gradient estimate assumption (\ref{C2_assumption}) in the $C^2$ estimate.
%
%\begin{theorem}[Fu-Yau \cite{FY}] \label{FY_estimate}
%Let the dimension of $X$ be equal to $n=2$. For any $\delta >0$, there exists $A_\delta>0$ depending on $(X,\o)$, $f$ $\alpha$, $\mu$, such that if $A < A_\delta$, then for any solution $u$ of the Fu-Yau equation (\ref{scalar_eqn}) with normalized condition (\ref{A_normalize}), there holds
%\be
%e^{-u} |Du|^2 \leq \delta.
%\ee
%\end{theorem} 
%{\it Proof.} 
From the previous analysis, we may choose $A_\delta$ such that 
\be
\kappa={\rm min}_X(e^{-2u}F) \geq 1 - \alpha \, \delta.
\ee
From (\ref{first_estimate}), we have
\be
1- \alpha \, \delta \leq e^{-2u} F \leq 1 - 4 \alpha e^{-u} |D u|^2 \bigg \{ 1 - (\|f\|_\infty +1) e^{-2u} \bigg \} + O(e^{-2u}).
\ee
After choosing to be $A_\delta$ smaller if necessary, we see that the previous inequality implies
\be \label{sharp_grad}
e^{-u} |Du|^2 \leq \delta.
\ee

\subsection{The case of higher dimension $n$}

%In higher dimensions, it is not difficult to see that the inequality from Theorem \ref{minimum} is not powerful enough to provide a lower bound for $\kappa_p$. In fact, even if we restrict ourselves only to the leading terms, obtained by setting formally $\theta=0$, it is easy to construct sequences of values $\kappa_p$ which satisfy the condition that the combination of leading terms be non-negative, but which tend to $0$.

In higher dimensions, it is not difficult to see that the inequality (\ref{minimum_ineq}) obtained in Theorem \ref{minimum} is not powerful enough to provide a lower bound for $\kappa_p$. In fact, even if we restrict ourselves only to the leading terms by setting formally $\theta=0$, the computation and examples indicate that the case $n=2$ case is quite special. In the $n=2$ case, as shown in (\ref{inequalitykp}) with $\theta=0$, it is easy to see that the leading terms about $\kappa_p$ cancel perfectly between both sides. However, this is not the case for higher dimensions.
\par
We illustrate the problem in the case $n=3$. Suppose that at the point $p \in X$ where $G$ achieves its minimum, $Du$ happens to be in the direction of $\lambda'_1$. Substituting $\kappa_c =3$ and $n=3$, the inequality (\ref{minimum_ineq}) with $\theta=0$ obtained in Theorem \ref{minimum} becomes
\bea
0 &\leq& \bigg( {1 \over 9} - {\kappa_p \over 3} \bigg) ( \lambda_1' + \lambda_2' + \lambda_3') e^{-u} + \bigg( {2 \over 9} \kappa_p - {4 \over 9} \bigg) e^{-u} ( \lambda_2' + \lambda_3') + {2 \over 3} \kappa_p \nonumber\\
&&+ {4 \over 9} e^{-3u} \lambda'_1 \lambda'_2 \lambda'_3 + \bigg({\kappa_p \over 9} - {5 \over 9}  \bigg) e^{-2u} \lambda'_2 \lambda'_3. 
\eea
This inequality cannot prevent $\kappa_p = e^{-2u} \sigma_2(\lambda')$ from starting at $\kappa_c=3$ and then going to zero along the method of continuity. Indeed, the path $e^{-u} \lambda'=(1, s, s)$ gives $\kappa_p = 2s + s^2$ and the previous inequality reduces to
\be
0 \leq {1 \over 9} (s^4 - 2s^2 +1) .
\ee
Thus it is unclear whether the non-degeneracy estimate holds in higher dimensions, and it would certainly require a different method.
%\par
%Roughly speaking, we can consider the inequality (\ref{inequality}) and try and use a similar argument as in the inequality (\ref{gradient1}) in order to bound $\left|Du\right|_F^2$ and $\sum_j \sigma_2(\lambda' | j) \frac{|u_j|^2}{|Du|^2}$. This will give us an inequality involving $\lambda_1', \sigma_1(\lambda'), \sigma_3(\lambda')$ and $\kappa_p$. Then, we can apply the standard Newton-MacLaurin inequality and Garding inequality to obtain an inequality about $\kappa_p$ which is analogous to (\ref{inequalitykp}), but more complicated. In the $n=2$ case, as shown in (\ref{inequalitykp}) with $\theta=0$, it is easy to see that the leading terms about $\kappa_p$ cancel perfectly between both sides. However, this is not the case for higher dimensions. Indeed, we can construct sequences of values $\kappa_p$ which satisfy the analogous inequality of (\ref{inequalitykp}) even for $n=3$, but tend to $0$. Thus it is unclear whether the non-degeneracy estimate holds in higher dimensions, and it would certainly require a different method.

\bigskip

\noindent {\bf Acknowledgements:} The authors would like to thank Pengfei Guan for stimulating conversations and for his notes on Fu-Yau's equation. They would also like to thank Valentino Tosatti for his lectures and notes on Strominger systems. The authors are also very grateful to the referee for a particularly careful reading of the paper, and for numerous suggestions which helped clarify the paper a great deal.

\bigskip

Department of Mathematics, Columbia University, New York, NY 10027, USA

\smallskip

phong@math.columbia.edu

\smallskip
Department of Mathematics, Columbia University, New York, NY 10027, USA

\smallskip
 picard@math.columbia.edu

\smallskip
Department of Mathematics, University of California, Irvine, CA 92697, USA

\smallskip
xiangwen@math.uci.edu

\end{document}